\documentclass[10pt]{article} 
\usepackage{a4,amssymb, amsmath, amsthm}
\usepackage{graphics}

\usepackage{subfig}
\usepackage{graphicx}
\usepackage{stackengine}

\newtheorem{theorem}{Theorem}[section]
\newtheorem{corollary}[theorem]{Corollary}
\newtheorem{lemma}[theorem]{Lemma}

\theoremstyle{definition}
\newtheorem{definition}[theorem]{Definition}

\theoremstyle{rem}

\newcommand{\basisphi}{\varphi_j^p}
\newcommand{\basisnu}{\nu_i^q}

\newcommand{\basisphilin}{\varphi_j}
\newcommand{\basisnulin}{\nu_i}

\newcommand{\ansatzVx}{V_h^p}
\newcommand{\ansatzVt}{V^q_{\Delta t}}
\newcommand{\ansatzVxt}{V^{p,q}}
\newcommand{\ansatzVxtgeneral}{V^{p,q}}

\def\XXint#1#2#3{{\setbox0=\hbox{$#1{#2#3}{\int}$ }
\vcenter{\hbox{$#2#3$ }}\kern-.6\wd0}}

\begin{document}

\providecommand{\keywords}[1]{{\noindent \textit{Key words:}} #1}
\providecommand{\msc}[1]{{\noindent \textit{Mathematics Subject Classification:}} #1}

\title{Adaptive time-domain boundary element methods for the  wave equation with Neumann boundary conditions}
\author{*°A. Aimi and *°G. Di Credico and $^\S$H. Gimperlein and *°C. Guardasoni\\\\\small *Dept.~of Mathematical Physical and Computer Sciences, University of Parma, Italy\\\small °Members of the INdAM-GNCS Research Group, Italy\\\small$^\S$Engineering Mathematics, University of Innsbruck, Austria}
\date{\today}

\maketitle \vskip 0.5cm
\begin{abstract}
\noindent This article investigates adaptive mesh refinement procedures for the time-domain wave equation with Neumann boundary conditions, formulated as an  equivalent hypersingular boundary integral equation. Space-adaptive and time-adaptive versions of a space-time boundary element method are presented, based on a reliable a posteriori error estimate of residual type. Numerical experiments illustrate the performance of the proposed approach.
\end{abstract}

\section{Introduction}

For elliptic problems,  adaptive versions of both finite element and boundary element methods give rise to fast approximations to nonsmooth solutions \cite{bonito,gwinsteph}.  Based on local error indicators obtained from an a posteriori error estimate, the adaptive algorithm produces a sequence of locally refined meshes. In many situations the resulting approximations were shown to have optimal convergence rates. For space-time boundary element formulations of time-dependent problems adaptive methods have been of much recent interest \cite{gantner,review,adaptive,Glaefke,hoonhout2023,zank1d}, building on the advances for time independent problems \cite{gwinsteph,feischl2015,cmsp,cms,cs}.\\
In this article we investigate adaptive mesh refinement procedures for the Neumann boundary value problem for the wave equation, formulated as a boundary integral equation in the time-domain. We present a reliable a posteriori error estimate of residual type for a large class of conforming discretizations. 
The error estimate defines local error indicators, which allow to define  both space-adaptive mesh refinement algorithms and the {adaptive selection of time steps}. For scattering problems from an open arc, the space-adaptive algorithm recovers the convergence rates known for time-independent problems.  \\

To describe the main results, we consider the wave equation
\begin{equation}\label{ibvp1}\partial_t^2 u - \Delta u = 0\ , \quad u=0\ \text{ for }\ t\leq 0\ ,\end{equation}
in the complement $\mathbb{R}^d \setminus \overline{\Omega}$ of a polyhedral domain or screen, $d=2,3$.
On the boundary $\Gamma = \partial \Omega$ we consider inhomogeneous Neumann boundary conditions,
\begin{equation}\label{ibvp2}\frac{\partial u}{\partial n} = f \ ,\end{equation}
with a given $f$ and the outer unit normal vector $n$ to $\Gamma$. To reduce this problem to the boundary, we use a double layer potential ansatz for $u$, 
\begin{align*}
u(t,x)&=\iint_{\mathbb{R}^+ \times \Gamma}  \frac{\partial G}{\partial n_y}(t- \tau,x,y)\ \psi(\tau,y)\, ds_y\, d\tau
\end{align*}
in terms of the fundamental solution of the wave equation, \begin{align}\label{eq:green_intro}
 G(t-\tau,x,y)&= \frac{H(t-\tau-|x-y|)}{2\pi \sqrt{(t-\tau)^2+|x-y|^2}},\qquad &\; d=2,\\
 G(t-\tau,x,y)&=\frac{\delta(t-\tau-|x-y|)}{4\pi |x-y|}, & d=3.
\end{align}
The density $\psi$ on $\mathbb{R}^+ \times \Gamma$ here satisfies $\psi(\tau,y) = 0$ for $\tau\leq 0$, $H$ denotes the Heaviside function and $\delta$ the Dirac distribution. The Neumann problem \eqref{ibvp1}, \eqref{ibvp2} is then equivalent to a time dependent hypersingular integral equation for $\psi$ :
\begin{align}\label{hypersingeq_intro}
\mathcal{W} \psi(t,x) = \iint_{\mathbb{R}^+\times \Gamma} \frac{\partial^2 G}{\partial n_x \partial n_y}(t- \tau,x,y)\ \psi(\tau,y)\ ds_y \ d\tau = f(t,x)\,.
\end{align}
The space-time boundary integral operator defined as in \eqref{hypersingeq_intro} arises from double layer ansatz for the solution to the wave equation, as stated in \cite{bh} for 3D problems. A similar definition can be found also in \cite{chen2010boundary} for the time independent case.\\
We consider space-time Galerkin discretizations of \eqref{hypersingeq_intro} in subspaces $V$ of piecewise polynomials in space and time, defined in Section \ref{sec:discretization}. We present a residual error estimate for the resulting space-time Galerkin discretizations,
which leads to  adaptive mesh refinement procedures, based on the four steps:
 \begin{align*}
  \textbf{SOLVE}&\longrightarrow \textbf{ESTIMATE}\longrightarrow \textbf{MARK}\longrightarrow \textbf{REFINE}.
\end{align*}
The numerical experiments study both space-adaptive and time-adaptive procedures for $d=2$ and they confirm the efficiency and reliability of the estimate. 
In particular, they include screen problems, where the geometric singularities pose the greatest numerical challenges, as well as transient singular behavior. \\

Our work builds on \cite{adaptive}  for the Dirichlet problem for the wave equation, formulated as a weakly singular integral equation, where an a posteriori estimate and a space-adaptive algorithm were considered. We also refer to \cite{Glaefke} for adaptive numerical experiments. As in \cite{adaptive}, our analysis relies on suitable weighted Sobolev spaces in time. Closely related function spaces were introduced also for finite element discretizations of the wave equation, assuming sufficient regularity in time, see \cite{cf23, chaumont2024damped} for the a posteriori error analysis of the Dirichlet problem. We are not aware of any related works for the Neumann problem or its formulation as a time dependent hypersingular integral equation. \\
The space-adaptive methods we introduce for the wave equation are particularly relevant for the efficient approximation of time-independent features, including singularities of the solution near edges and corners. For such geometric singularities, quasi-optimal convergence rates on time-independent graded meshes were obtained using finite elements in \cite{mueller} for polygons in $\mathbb{R}^2$ and using boundary elements in \cite{graded,hp} for edge and cone singularities in $\mathbb{R}^3$. Space-adaptive boundary element methods for the Dirichlet problem in these geometries, with a uniform time step, were studied in \cite{adaptive}. There the improved convergence rates known for adaptive methods for the Laplace equation were recovered. First
computational results towards fully space-time refinements were presented by Gl\"{a}fke \cite{Glaefke} for the Dirichlet problem in $\mathbb{R}^2$, and unpublished work
by Abboud uses ZZ error indicators towards space-adaptive mesh refinements
for screen problems in $\mathbb{R}^3$. \\

Previous work on adaptivity in time for boundary element methods is more limited. For the Dirichlet problem in $\mathbb{R}$ the adaptive selection of time steps to resolve singular temporal behavior was studied in \cite{hoonhout2023,zank1d}, following a first work \cite{sv} in $\mathbb{R}^3$. Also for time discretizations using convolution quadrature, instead of a Galerkin discretization, non-adaptively chosen, variable time steps have recently been of much interest, following \cite{ls}. \\
The techniques used here build on the numerical analysis of adaptive boundary element methods for the Laplace equation. We particularly refer to \cite{cmsp,cs} for the hypersingular integral equation and to \cite{gwinsteph} for an overview. \\
Our work contributes to the recent wider interests in boundary element methods for wave equations, see \cite{costabel2004time,hd,review,sayas2016} for an overview of both Galerkin and convolution quadrature methods. Such time domain methods are of particular relevance for problems which cannot be reduced to the frequency domain, including nonlinear problems and problems that involve a broad range of frequencies \cite{aimi2023time,aimi2024space,contact}. We note recent, complementary works on the efficient assembly and compression of the space-time matrices for both time-stepping and more general space-time discretizations \cite{aimi2020,hsiao1,bertoluzza,merta}.
The numerical analysis of the Neumann problem in this paper goes back to \cite{bh1986}. We refer to \cite{becache1993variational,becache1994space} for early investigations involving geometric singularities on screens and to \cite{banz,nedelec} for recent applications.\\

\noindent \emph{Structure of this article:} Section \ref{sec:problem} recalls the formulation of the Neumann problem for the wave equation as a time-dependent hypersingular integral equation and discusses its weak formulation.  In Section \ref{sec:discretization} the weak formulation is discretized in space and time using tensor products of piecewise polynomials. Section \ref{sec:apost} derives the a posteriori error estimate,
while Section \ref{algo} discusses the algorithmic implementation of the space-time system of linear equations and of the residual error indicators. The adaptive algorithms are introduced and studied in Section \ref{NE}: Subsection \ref{subsection:spaceadaptive} presents a space-adaptive mesh refinement procedure followed by several numerical experiments; finally, a procedure which adaptively selects the time steps is introduced and numerically studied in Subsection \ref{subsection:timeadaptive}. \\

\noindent \emph{Notation:} We write $f \lesssim g$ provided there exists a constant $C$ such that $f \leq Cg$. If the constant
$C$ is allowed to depend on a parameter $\sigma$, we write $f \lesssim_\sigma g$.

\section{Problem formulation}\label{sec:problem}

We recall from the Introduction the formulation of the Neumann problem \eqref{ibvp1}, \eqref{ibvp2} for the wave equation as a time-dependent hypersingular boundary integral equation 
\begin{equation}\label{hypersingeq}
 \mathcal{W} \psi = f\ .
\end{equation}
Its weak formulation involves the bilinear form \begin{equation}B(\psi, \phi) := \langle \mathcal{W} {\partial_t} \psi,\phi \rangle_\sigma := \iint_{\mathbb{R}^+ \times \Gamma} \mathcal{W} \partial_t {\psi}(t,x)\ \phi(t,x)\ ds_x \, d_\sigma t,\end{equation} where $d_\sigma t = \mathrm{e}^{-2\sigma t} dt$ for fixed $\sigma>0$. The properties of $B(\cdot,\cdot)$ follow from the mapping and coercivity properties of $\mathcal{W}$ in the space-time anisotropic Sobolev spaces $H_\sigma^{r}(\mathbb{R}^+, H^{s}(\Gamma))$, whose definition is recalled in the Appendix:  
\begin{theorem}\label{mapthm}
a) Let $r \in \mathbb{R}$. Then the hypersingular operator $\mathcal{W}  $ is continuous 
\begin{align*}
\mathcal{W}  &: H_\sigma^{r+1}(\mathbb{R}^+, \widetilde{H}^{\frac{1}{2}}(\Gamma)) \to H_\sigma^{r}(\mathbb{R}^+, H^{-\frac{1}{2}}(\Gamma))\ .
\end{align*}
b) The operator $\mathcal{W}  \partial_t$ is weakly coercive: $$\iint_{\mathbb{R}^+ \times \Gamma}(\mathcal{W} {\partial_t} \psi(t,x))\,  \psi(t,x)\  ds_x \ d_\sigma t \gtrsim_\sigma \|\psi\|_{0,{\frac{1}{2}},\Gamma,\ast}^2.$$ 
It follows that $\mathcal{W}^{-1} : H_\sigma^{r+1}(\mathbb{R}^+, \widetilde{H}^{-\frac{1}{2}}(\Gamma)) \to H_\sigma^{r}(\mathbb{R}^+, H^{\frac{1}{2}}(\Gamma))$ is continuous.  
\end{theorem}
\noindent For the proof, see \cite{hd} for part a) when $\partial \Gamma=\emptyset$. In this case part b) follows from equation (2.14), p.~174 in \cite{thd}. For the half-space or when $\partial \Gamma \neq\emptyset$, a) is shown in \cite{setup}; the proof of b) is obtained by extending Ha Duong's proof in \cite{hd} for $\partial \Gamma=\emptyset$, using the modifications from \cite{setup}. \\

We conclude that $B(\cdot,\cdot)$ is continuous and weakly coercive: \begin{corollary}\label{NPbounds} For every $\phi,\psi \in H^1_\sigma( \mathbb{R}^+, \widetilde{H}^{\frac{1}{2}}(\Gamma))$ there holds:
$$|B(\psi,\phi)| \lesssim \|\psi\|_{1,\frac{1}{2},\Gamma, \ast} \|\phi\|_{1,\frac{1}{2}, \Gamma,\ast}$$
and
$$\|\psi\|_{0,\frac{1}{2},\Gamma,\ast}^2 \lesssim B(\psi,\psi) . $$
\end{corollary}
Note the loss of a time derivative between the upper and lower estimates. We refer to \cite{jr} for a detailed investigation of mapping and coercivity properties and to \cite{sut} for an investigation of an alternative inf-sup stable bilinear form.  \\  

Using the bilinear form $B(\cdot,\cdot)$, we obtain the following weak formulation of equation \eqref{hypersingeq}: \\

\noindent \textit{find} $\psi \in H^1_\sigma( \mathbb{R}^+, \widetilde{H}^{\frac{1}{2}}(\Gamma))$ \textit{such that} 
\begin{equation}\label{NP}
B(\psi,\phi)=\langle \partial_t f,\phi\rangle_\sigma\, ,\quad \text{ for all }\phi \in H^1_\sigma( \mathbb{R}^+, H^{\frac{1}{2}}(\Gamma)).
\end{equation}

Then, we consider conforming Galerkin discretizations in a subspace $V \subset H^1_\sigma( \mathbb{R}^+, \widetilde{H}^{\frac{1}{2}}(\Gamma))$:\\

\noindent \textit{find} $\psi_{h,\Delta t} \in V$ \textit{such that} 
\begin{equation}\label{NPdisc}
B(\psi_{h,\Delta t},\phi_{h,\Delta t})=\langle \partial_t f,\phi_{h,\Delta t}\rangle_\sigma\, ,\quad \text{ for all } \phi_{h, \Delta t} \in V.
\end{equation}

The well-posedness of the continuous and discretized problems are a basic consequence of Corollary \ref{NPbounds}:
\begin{corollary}
Let $f \in H^2_\sigma(\mathbb{R}^+,H^{-\frac{1}{2}}(\Gamma))$. Then the Neumann problem \eqref{NP} and its discretization \eqref{NPdisc} admit unique solutions $\psi \in H^1_\sigma(\mathbb{R}^+,\widetilde{H}^{\frac{1}{2}}(\Gamma))$, $\psi_{h,\Delta t} \in V$, and the estimates $$\|\psi\|_{1,\frac{1}{2},\Gamma,\ast}\lesssim \|f\|_{2,-\frac{1}{2},\Gamma},\quad \|\psi_{h,\Delta t}\|_{1,\frac{1}{2},\Gamma,\ast} \lesssim \|f\|_{2,-\frac{1}{2},\Gamma},$$ hold.
\end{corollary}

\section{Discretization}\label{sec:discretization}

For simplicity we assume that the boundary $\Gamma$ is {polygonal for $d=2$ or polyhedral for $d=3$. 
We introduce a polyhedral space-time mesh ${\cal T}$ on $\mathbb{R}^+ \times \Gamma$ consisting of elements $\tau_j$ (triangles or quadrilaterals for $d=2$, tetrahedra or prismatic elements for $d=3$), such that $\mathbb{R}^+ \times\Gamma=\cup_{j}^{N} \tau_j$. Here, we assume that $\tau_i \cap \tau_j=\emptyset$ if $i \neq j$. For $\tau \in \mathcal{T}$ we denote by $H_\tau$ the diameter of $\tau$. The diameters $H_\tau$ define a function $H_{\mathcal{T}}$ by $H_{\mathcal{T}}(x,t) = H_\tau$ when $(x,t) \in \tau$. 
On the mesh ${\cal T}$ we consider the space $\ansatzVxtgeneral \subset H^1_\sigma( \mathbb{R}^+, \widetilde{H}^{\frac{1}{2}}(\Gamma))$ of piecewise polynomial continuous functions of a fixed degree $p{\geq 1}$ in space, of a fixed degree $q{\geq 1}$ in time. A basis for the space $\ansatzVxtgeneral$ is denoted by $\{\varphi_j^{p,q}\}$.The Galerkin discretization \eqref{NPdisc} is considered in the spaces $V=\ansatzVxt$.\\
Different degrees $p, q$ are particularly relevant in the special case of tensor product discretizations, 
where the space-time domain $\mathbb{R}^+ \times \Gamma$ is discretized separately in space and in time. As the numerical experiments in this article particularly consider such space-time meshes, we describe them in more detail.\\
We introduce a mesh ${\cal T}_\Gamma$ on $\Gamma$ consisting of elements $\Gamma_j$ (straight segments for $d=2$, triangular faces for $d=3$), such that $\Gamma=\cup_{i=1}^{N_{\Gamma}} \Gamma_j$. We assume that $\Gamma_i \cap \Gamma_j=\emptyset$ if $i \neq j$. Denote by $h_j$ the diameter of $\Gamma_j$ and let $h=\max_j h_j$. On the mesh ${\cal T}_\Gamma$ we consider the space $\ansatzVx$ of piecewise polynomial continuous functions of a fixed degree $p{\geq 1}$. A basis for the space $\ansatzVx$ is denoted by $\{\basisphi\}$.  \\
For the time discretization we consider a decomposition ${\cal T}_{\mathbb{R}^+}$ of the time interval $\mathbb{R}^+$ into subintervals $I_i=[t_{i-1}, t_i)$ with time step $|I_i|=\Delta t_i$, $i\geq 1$. Let $\Delta t = \sup_i \Delta t_i$. On the mesh ${\cal T}_{\mathbb{R}^+}$ we consider the space $\ansatzVt$ of piecewise  polynomial  continuous functions, vanishing at $t=0$, of degree of $q\geq 1$. A basis for the space $\ansatzVt$ is denoted by $\{\basisnu\}$.

From ${\cal T}_{\mathbb{R}^+}$ and ${\cal T}_\Gamma$, we obtain on $\mathbb{R}^+ \times \Gamma$ a space-time mesh $\mathcal{T} = {\cal T}_{\mathbb{R}^+} \times {\cal T}_\Gamma =\bigcup_{i,j} \Box_{i,j}$, where $\Box_{i,j} = I_i \times  \Gamma_j$.
The diameter of a space-time element $\tau=\Box_{i,j}$ is of size $H_\tau \simeq \sqrt{\Delta t_i^2 + h_j^2} \simeq \max\{\Delta t_i, h_j\}$. Note that the temporal mesh is the same for all elements in space, respectively, the spatial mesh is the same for all time intervals. 
A basis for $\ansatzVxt =\ansatzVx \otimes \ansatzVt$  is given by tensor products $\basisnu(t) \basisphi(x)$. The Galerkin approximation $\psi_{h, \Delta t}$ of the solution $\psi$ hence takes the form
\begin{equation}\label{discrete psi sec3}
\psi_{h, \Delta t}(x,t):=\sum_{i} \sum_{j}
\alpha_{j}^{(i)} \,\basisphi(x)\,\basisnu(t).
\end{equation}

\section{A posteriori error estimate}\label{sec:apost}

{We adapt a construction originally due to \cite{cb}. Denote by $\mathcal{N}$ the set of all nodes in the space-time mesh $\mathcal{T}$. The nodes which lie in the interior of $\mathbb{R}^+\times \Gamma$ are denoted by $\mathcal{K}$. For each $z\in \mathcal{K}$ we write $\varphi_z \in \ansatzVxt$  for a basis function with $\varphi_z(z) = 1$, which equals $0$ in all nodes $x \in \mathcal{N} \setminus \{z\}$. }

{From the basis functions $\varphi_z$, for $p=1$ we construct a partition of unity   $\psi_z$, i.e., satisfying $\sum_{z\in \mathcal{K}} \psi_z = 1$.  To do so, we define $\zeta : \mathcal{N} \to \mathcal{K}$ as follows: For $z \in \mathcal{K}$, we set $\zeta(z) = z$, while for $z \in \mathcal{N} \setminus \mathcal{K}$ we choose any node $\zeta(z) \in \mathcal{K}$. Denoting the preimage of $z \in \mathcal{K}$  by $I(z)$, $I(z) = \{\tilde{z} : \zeta(\tilde{z}) = z\}$, we define
$$\psi_z = \sum_{\tilde{z} \in I(z)} \phi_{\tilde{z}}.$$ One readily checks that $\sum_{z\in \mathcal{K}}\psi_z = 1$.
The interior of the support of $\psi_z$ is denoted by $\Gamma_z$: $\Gamma_z = \{(x,t) \in \mathbb{R}^+\times \Gamma : \psi_z(x,t)>0\}$. }\\

\noindent {\textbf{Assumption:} $\Gamma_z$ is connected and a union of finitely many elements $\tau \in \mathcal{T}$.}\\

\noindent We define $H_z$ as the diameter of the largest element in $\Gamma_z$: $H_z = \max\{H_\tau : \tau \in \Gamma_z\}$.\\
Proceeding as in the time-independent case, \cite{cb}, Section 2, we obtain:
\begin{lemma}\label{cbeq210}
There exists a constant $C$ that depends on $\Gamma$ and the aspect ratio of the space-time elements (but not their sizes) such that for all $z \in \mathcal{K}$ and $g \in H^{r}_\sigma(\mathbb{R}^+,\widetilde{H}^{1}(\Gamma))$: 
$$\|\psi_z g - \varphi_z g\|_{r,0,\Gamma_z} \leq C \min\{\|g\|_{r,0,\Gamma_z}, H_z\|g\|_{r,1,\Gamma_z} \}.$$
\end{lemma}

{We also note the following observation from \cite{cms}, Lemma 3.1, as adapted to time-dependent problems in \cite{adaptive}, Lemma 4.3.
\begin{lemma}\label{cc31}
Let $\Phi = \{\phi_j\}_j$, $\sum_{j} \phi_j = 1$, be a locally finite partition of unity of $\Gamma$ with overlap $K(\Phi) = \max_{j} \mathrm{card}\{k : \phi_k \phi_j \neq 0\}<\infty$. Then there exists a partition of $\Phi$ into $K \leq K(\Phi)$ non-empty subsets $\Phi_1, \dots, \Phi_K$, such that $\bigcup_{j=1}^K \Phi_j = \Phi$, $\Phi_j \cap \Phi_k = \emptyset$ if $j \neq k$ and for all $l \in \{1,\dots, K\}$ and $\phi_j,\phi_k \in \Phi_l$ with $j \neq k$, $\phi_j \phi_k = 0$ on $\mathbb{R}^+\times \Gamma$.
\end{lemma}}

{From Lemma \ref{cbeq210} and Lemma \ref{cc31} we obtain as in \cite{gwinsteph}, Section 10.4, a preliminary a posteriori error estimate:
\begin{lemma}\label{lem:localize}
For any $\mathcal{R} \in H^{r+1}_\sigma(\mathbb{R}^+,{H}^{0}(\Gamma))$ with $\langle\partial_t \mathcal{R} ,\varphi_z \rangle_\sigma =0$ for all $z \in \mathcal{K}$  and all $s \in [0,1]$:
$$\|\mathcal{R}\|_{r+1,s-1, \Gamma} \lesssim \left(\sum_{z \in \mathcal{K}} H_z ^{2-2s} \|\mathcal{R}\|_{r+1,0, \Gamma_z}^2  \right)^{1/2} \lesssim \|H_{\mathcal{T}}^{1-s}\mathcal{R}\|_{r+1,0, \Gamma}.$$
\end{lemma}}

{We will apply Lemma \ref{lem:localize} to the residual $\mathcal{R}:=f-\mathcal{W} \psi_{h, \Delta t}$. Using $\psi-\psi_{h, \Delta t} = \mathcal{W}^{-1} (f-\mathcal{W}\psi_{h, \Delta t}) = \mathcal{W}^{-1} \mathcal{R}$ and the boundedness of $\mathcal{W}^{-1}$ from Theorem \ref{mapthm}, we find
$$
\|\psi-\psi_{h, \Delta t}\|_{r, s, \Gamma, \ast}  = \|\mathcal{W}^{-1} \mathcal{R}\|_{r, s, \Gamma, \ast} 
 \lesssim_\sigma \|\mathcal{R}\|_{r+1, s-1, \Gamma}.
$$
Lemma \ref{lem:localize} gives a bound for the right hand side, leading to the following a posteriori error estimate. 
\begin{theorem}\label{apostthm}  Fix $r \in [-1,1]$, $s\in [0,\frac{1}{2}]$. Let $\psi  \in H^{1}_\sigma(\mathbb{R}^+,H^{\frac{1}{2}}(\Gamma))$ be the solution to \eqref{NP}, and let $\psi_{h,\Delta t} \in H^{1}_\sigma(\mathbb{R}^+,H^{\frac{1}{2}}(\Gamma))$. Then 
\begin{equation}\label{aposteriori_general}
\|\psi-\psi_{h, \Delta t}\|_{r, s, \Gamma, \ast} \lesssim_\sigma \left(\sum H_{z}^{2-2s} \|\mathcal{R}\|_{r+1,0, \Gamma_z}^2  \right)^{1/2} \lesssim \|H_{\mathcal{T}}^{1-s}\mathcal{R}\|_{r+1,0, \Gamma_z}.
\end{equation}
In the special case of tensor-product meshes and $r=-1$ we conclude
\begin{equation}\label{aposteriori}
\|\psi-\psi_{h, \Delta t}\|_{-1, s, \Gamma, \ast} ^2\lesssim_\sigma \sum_{\Box_{i,j} 
} 
\max\{\Delta t_i, h_j\}^{2-2s} \|\mathcal{R}\|_{0,0,\Box_{i,j}}^2.
\end{equation}
\end{theorem}
}
\noindent     Note that the constants in Theorem \ref{apostthm} depend on the aspect ratio of the space-time
elements, but not their sizes.\\

\noindent {\bf Remark.} 
Introducing the local error indicator as
\begin{equation}\label{local error ind}
    \eta(\Box_{i,j})=\max\{\Delta t_i, h_j\}^{2-2s} \|\mathcal{R}\|_{0,0,\Box_{i,j}}^2,
\end{equation}
we recall from Section \ref{sec:discretization} that in \eqref{local error ind} $\max\{\Delta t_i,h_j\}$ can be replaced by $\sqrt{\Delta t_i^2+h_j^2}$. In fact, the resulting estimate is equivalent, because  \eqref{aposteriori} contains a constant factor hidden in the notation $\lesssim$ and because $\max\{\Delta t_i,h_j\}\leq \sqrt{\Delta t_i^2+h_j^2}\leq\sqrt{2}\max\{\Delta t_i,h_j\}$, i.e., the two quantities differ at most by a constant factor $\sqrt{2}$.\\ 

\noindent {\bf Remark.} {The right hand side of formula \eqref{aposteriori} furnishes the error indicator of the adaptive algorithms described in Section \ref{NE}. In particular, this requires the evaluation on each space-time box $\Box_{i,j} = [t_{i-1},t_i) \times  \Gamma_j$ of the integral of the squared residual error 
$$\|\mathcal{R}\|_{0,0,  {\Box_{i,j}}}^2=\int_{t_{i-1}}^{t_i}\int_{\Gamma_j}[f(t,x)-\mathcal{W} \psi_{h, \Delta t}(t,x)]^2 ds_x \, d_\sigma t,$$
which can be managed by a suitable space-time quadrature formula, once fixed $\sigma=0$.
The computation of the squared residual error rests on the numerical approximation of $\mathcal{W} \psi_{h, \Delta t}(t,x)$, where $(t,x)$ will be each of the quadrature nodes related to the chosen numerical integration rule. Details of this process are described in Subsection \ref{algo_sub2}.\\
Let us finally note that, as in all adaptive methods, the evaluation of the local error indicators represents an extra cost with respect to a uniform space-time refinement procedure. Nevertheless, the adaptive algorithm is expected to attain the desired accuracy at lower numbers of DoFs and consequently with lower memory requirements and CPU times. An analysis of memory savings and examples of run times of the adaptive algorithm will be given in Section \ref{NE}.

\section{Implementation details} \label{algo}

In view of the numerical results given in Section \ref{NE}, here we  focus on the algebraic formulation of  the Galerkin discretization \eqref{NPdisc} with an emphasis on $d=2$. To simplify the exposition, we describe the case of tensor product discretizations, where a solution to \eqref{NPdisc} is sought in the space $\ansatzVxt= \ansatzVx \otimes \ansatzVt$. 

\subsection{Structure of space-time matrix}\label{algo_sub}

\noindent The details of the implementation in this setting are explained for uniform time steps in \cite{Aimi2010}, and we now discuss their generalization.
Specifying the discretization from Section \ref{sec:discretization}, on the bounded time interval $[0,T] = \bigcup_{i=0}^{N_T-1} [t_i, t_{i+1})$ we have $M_T=qN_T$
piecewise polynomial basis functions $\basisnu(t)$. Similarly, there are $M_\Gamma$ piecewise polynomial basis functions $\basisphi(x)$ on the mesh $\Gamma = \bigcup_{j=1}^{N_\Gamma}\Gamma_j$.
The Galerkin approximation $\psi_{h, \Delta t}$ of the solution $\psi$ then takes the form\\
\begin{equation}\label{discrete psi}
\psi_{h, \Delta t}(x,t):=\sum_{i=0}^{M_T-1} \sum_{j=1}^{M_\Gamma}
\alpha_j^{(i)} \,\basisphi(x)\,\basisnu(t)
\end{equation}
and the  discretized weak equation \eqref{NPdisc} leads to the linear system
\begin{equation}\label{linear_system}
\mathbb{E}\,\pmb{\alpha}=\pmb{\beta}.
\end{equation}
The matrix $\mathbb{E}$ has a time block lower triangular 
structure, since its elements vanish if $t_{\tilde{i}}< t_{i}$. Each time block $\mathbb{E}_{\tilde{i},\,i}$ has
dimension $qM_\Gamma$. The linear system can be written as
\begin{equation}\label{linear_system_explicit}
\begin{pmatrix}
\mathbb{E}_{0,0} & \cdots            & 0                 & \ldots    & 0\\
\vdots & \ddots            & 0                 & \ldots    & 0\\
\mathbb{E}_{i,0} & \cdots & \mathbb{E}_{i,i}                & \ldots    & 0\\
\vdots        & \ddots          & \vdots            & \ddots    & \vdots\\
\mathbb{E}_{N_T-1,0}   & \cdots  & \mathbb{E}_{N_T-1,i}   & \ldots    & \mathbb{E}_{N_T-1,N_T-1} 
\end{pmatrix}
\begin{pmatrix}
\pmb{\alpha}^{(0)}\\
\vdots  \\
\pmb{\alpha}^{(i)}\\
\vdots       \\
\pmb{\alpha}^{(N_T-1)}
\end{pmatrix}=
\begin{pmatrix}
\pmb{\beta}^{(0)}\\
\vdots \\
\pmb{\beta}^{(i)}\\
\vdots       \\
\pmb{\beta}^{(N_T-1)}
\end{pmatrix}
\end{equation}
where the unknowns and rhs entries are reorganized as follows, for $i=0,\cdots,N_T-1$
\begin{equation}\label{alpha_beta}
\pmb{\alpha}^{(i)}=\left(\alpha_{i ,1},\cdots,\alpha_{i,
qM_\Gamma}\right)^\top, \;\pmb{\beta}^{(i)}=\left(\beta_{i,
1},\cdots,\beta_{i, qM_\Gamma}\right)^\top\,. 
\end{equation}
The solution of $(\ref{linear_system_explicit})$ is obtained by a
block forward substitution, since $\mathbb{E}$ has non singular diagonal blocks.\\

\noindent {\bf Remark.} When the space-time mesh is a global product of spatial and temporal meshes with \emph{equal} time steps, the Galerkin matrix in \eqref{linear_system_explicit} has a block Toeplitz structure with constant blocks $\mathbb{E}_{\ell,0} = \mathbb{E}_{i,i'}$ corresponding to global time steps $t_i \geq t_{i'}$ with $\ell=i-i'$. The Toeplitz structure is lost when the time step or spatial mesh change with time. To minimize this computational overhead, efficient implementations of time-adaptive algorithms need to reuse those entries of the space-time system which have not been affected by the last refinement step.\\

We now describe the evaluation of the matrix entries in \eqref{linear_system_explicit} in the case of polynomial degrees $p=q=1$, which will be kept fixed in the remaining of the paper and in particular for the numerical experiments. With this choice, the piecewise linear temporal shape
functions, for $i=0,\cdots,N_T-1$, take the form
$$
\basisnulin(t)=R(t-t_i)-\,R(t-t_{i+1}),
$$
where $R(t-t_k)=\frac{t-t_k}{\Delta\,t}\,H[t-t_k]$ is a ramp
function.
Similarly, the piecewise linear hat functions on $\Gamma = \bigcup_{j=1}^{N_\Gamma}\Gamma_j$ are denoted by $\basisphilin(x)$.
Having set $\Delta_{\tilde{i}i}=t_{\tilde{i}}-t_i$, the elements of a generic matrix block, after a double analytic integration in the
time variables, are of the form 
\begin{equation}
e^{\tilde{i},i}_{\tilde{j},j}=\sum_{\gamma,\delta=0}^1(-1)^{\gamma+\delta}\!\!\int_{\Gamma}
\varphi_{\tilde{j}}(x)\!\!
\int_{\Gamma}H[\Delta_{\tilde{i}+\gamma,i+\delta}-|x-y|]\,
\mathcal{D}(x,y,t_{\tilde{i}+\gamma},t_{i+\delta})\,
\varphi_j(y)\,ds_y\, ds_x\,,
\label{App}
\end{equation}
where
\begin{eqnarray}
&&\mathcal{D}(x,y,t_{\tilde{i}},t_i)= \displaystyle
{\frac{1}{2\pi\,\Delta t_{\tilde{i}}\,\Delta t_i}\Big\{\frac{(x-y)\cdot
n_x\,\,\,(x-y) \cdot n_y}{r^2}\,\, \frac{\Delta_{\tilde{i}i}
\sqrt{\Delta_{\tilde{i}i}^2-|x-y|^2}}{|x-y|^2}\,+} \nonumber \\
&&\displaystyle{\frac{n_x \cdot
n_y}{2}
\Big[\log(\Delta_{\tilde{i}i}+\sqrt{\Delta_{\tilde{i}i}^2-|x-y|^2})-\log |x-y|
-\frac{\Delta_{\tilde{i}i}\,\sqrt{\Delta_{\tilde{i}i}^2-|x-y|^2}}{|x-y|^2}} \Big]\Big\}\,,
\label{calD}
\end{eqnarray}
and they are evaluated using the standard element by element technique and quadrature formulas explicitely developed for the wave equation hypersingular kernel in \cite{Aimi2010}.

\subsection{Implementation of error indicators}\label{algo_sub2}

Here we give practical indications for the implementation of the error indicator introduced in Section \ref{sec:apost}, which depends specifically on the computation of the residual $\mathcal{R}(x,t)=f(x,t)-\mathcal{W}\psi_{h,\Delta_t}(x,t)$. Recalling the definition of the hypersingular operator $\mathcal{W}$ in \eqref{hypersingeq_intro} and the definition of the discrete solution in \eqref{discrete psi} we write:
\begin{equation}\label{hyp in res}
\mathcal{W}\psi_{h,\Delta_t}(x,t)=\sum_{i=0}^{N_T-1} \sum_{j=1}^{M_\Gamma}\alpha_{j}^{(i)}\int_0^T\int_{\Gamma} \frac{\partial^2 G}{\partial n_x \partial n_y}(t- \tau,x,y)\  \,\varphi_j(y)\,\nu_i(\tau)\ ds_y\ d\tau\,.
\end{equation}
BEM discretization allows to reduce computation of \eqref{hyp in res} as a summation of space-time integrals depending of the single products of shape function $\varphi_j(y)\nu_i(\tau)$. Moreover, exploiting the definition of the ramp functions $\nu_i(\tau)$, analytical integration in time can be performed, leading to
\begin{equation}\label{hyp in res space 1}
\int_0^T\int_{\Gamma} \frac{\partial^2 G}{\partial n_x \partial n_y}(t- \tau,x,y)\  \,\varphi_j(y)\,\nu_i(\tau)\ ds_y\ d\tau =\frac{1}{2\pi \Delta t_i }\sum_{\delta=0}^1(-1)^\delta \int_{supp(\varphi_j)} \widetilde{\mathcal{D}}(x,y,t-t_{i+\delta})\varphi_j(y)\ ds_y
\end{equation}
where, indicating by $\tilde{\Delta}$ the difference between two generic time instants,
\begin{equation}\label{D_tilde}
\widetilde{\mathcal{D}}(x,y,\tilde{\Delta})=\frac{f(x,y,\tilde{\Delta})H[\tilde{\Delta}-\vert x-y \vert]}{\vert x-y \vert^2}\,,
\end{equation}
with regular part
$$f(x,y,\tilde{\Delta})=\left(-n_x \cdot n_y+\frac{(x-y)\cdot
n_x\,(x-y) \cdot n_y}{\vert x-y \vert^2}\right)\sqrt{\tilde{\Delta}^2-\vert x-y \vert^2}+\tilde{\Delta}^2\frac{(x-y)\cdot n_x\, (x-y) \cdot n_y}{\vert x-y \vert^2\sqrt{\tilde{\Delta}^2-\vert x-y \vert^2}}$$
and singularity $1/\vert x-y \vert^2$ of the hypersingular kernel highlighted. 
To numerically evaluate \eqref{hyp in res space 1} in the framework of the standard element-by-element technique, it is necessary to compute contributions of the type
\begin{equation}\label{fundamental spatial integral}
\int_{\Gamma_k}
\frac{f(x,y,\tilde{\Delta})H[\tilde{\Delta}-\vert x-y \vert]}{\vert x-y \vert^2}\varphi_j(y)\ ds_y
\end{equation}
where the space mesh element $\Gamma_k \subset supp(\varphi_j)$. 
If the field point $x \in \Gamma_k$ the integral is understood as a Hadamard finite part and managed by an appropriate quadrature formula as in \cite{{Aimi2010}}. Conversely, whenever the evaluation occurs for a field point $x$ outside $\Gamma_k$, a Gauss-Legendre rule is sufficient to achieve a good approximation. The discontinuity due to the Heaviside function in the kernel can be overcome splitting $\Gamma_k$ in additional subregions: partitions are determined by the discontinuities of the Heaviside function as seen in \cite{ADGMP2009}.

\section{Adaptive algorithms and numerical results}\label{NE}

In this Section, we will show several numerical results obtained applying the algorithms for adaptive space/time mesh refinements introduced below. 
The a posteriori error estimate from Theorem \ref{apostthm} leads to  adaptive procedures, based on the four steps:
 \begin{align*}
  \textbf{SOLVE}&\longrightarrow \textbf{ESTIMATE}\longrightarrow \textbf{MARK}\longrightarrow \textbf{REFINE}.
\end{align*}

\subsection{Space-adaptivity}\label{subsection:spaceadaptive}

As shown in \cite{aimi23,graded, mueller, mueller2}, for polyhedral domains and screens, time-independent graded meshes lead to quasi-optimal convergence rates in spite of the singular behavior of the solution. A corresponding space-adaptive boundary element procedure for the time domain weakly singular integral equation was considered in \cite{adaptive}, which adaptively refines the spatial mesh, but uses a uniform time step. \\
Analogously, as a step towards adaptive mesh refinements for the hypersingular integral equation in space and time, we first focus on refinements of the spatial mesh, based on time-integrated error indicators. The precise algorithm is given in the following:\\

\noindent {{\textbf{Space-Adaptive Algorithm} 
}}\\
 Input: Mesh $\mathcal{T}=(\mathcal{T}_T \times \mathcal{T}_\Gamma)_0$, refinement parameter $\theta  \in (0, 1)$, tolerance $\epsilon > 0$, datum $f$.
\begin{enumerate}
\item Solve \eqref{NPdisc} on $\mathcal{T}$.
\item Compute the error indicators $\eta(\Box_{i,j})$ 
in each space-time element $\Box_{i,j} \in \mathcal{T}$, i.e., for $i=0, \dots, N_T-1$ and $j=1,\dots, M_\Gamma$. 
\item Evaluate $\eta_j=\displaystyle \sum_{i=0}^{N_T-1} \eta(\Box_{i,j}), \, j=1,\dots,M_\Gamma$.
\item Stop if $\displaystyle \sum_{j=1,\dots,M_\Gamma} \eta_j < \epsilon$\,, otherwise:
\item Find $\displaystyle \eta_{max}=\max_{j=1,\dots,M_\Gamma} \,\eta_j$. 
\item Mark all $\Box_{i,j}\in \mathcal{T}$ with $j$ such that $\eta_j > \theta \,\eta_{max}$ and $i=0, \dots, N_T-1$. 
\item Refine each marked $\Box_{i,j}$ in space,  keeping $ \max_j \frac{\Delta t_i}{h_j}$ fixed.
\item Go to 1.
\end{enumerate}
Note that the space-adaptive algorithm preserves the maximal CFL ratio $\frac{\Delta t_i}{h_j}$. While the Galerkin discretization is stable without this constraint, it assures the reduction of the time discretization error. Also other choices for the marking are commonly used, e.g.~D\"{o}rfler marking. Instead of refining in time in Step 7., one can also consider space-adaptive algorithms for a sufficiently small, fixed time step.\\
For all the following numerical simulations, the parameter $\theta=0.4$ has been considered. 

\subsubsection{Space refinement step}
In the following we explain how the space-adaptive procedure influences the construction of new basis functions, focusing in particular on the generic subsequent 
basis functions $\varphi_{j-1}$ and $\varphi_{j}$ 
considered in Section \ref{algo_sub}, 
which share a common mesh element $\Gamma_j$ representing the intersection of their support. Let us suppose that the adaptive procedure compels its refinement: hence $\Gamma_j$ is split in two segments $\Gamma_j^{(1)}$ and $\Gamma_j^{(2)}$ and a new hat function $\tilde{\varphi}_{j-1,j}$ having $\Gamma_j=\Gamma_j^{(1)}\cup \Gamma_j^{(2)}$ as support is created. The new mesh elements lead also to the update of $\varphi_{j-1}$ and $\varphi_{j}$ which become the new hat functions $\tilde{\varphi}_{j-1}$ and $\tilde{\varphi}_{j}$, as visualized in Figure \ref{fig:spatial_mesh_ref_ex}.\\
The described procedure induces in all the temporal blocks $\mathbb{E}_{\tilde{i},i}$ the insertion of a new row $\tilde{e}_{(j-1,j),:}$ and a new column $\tilde{e}_{:,(j-1,j)}$  related to the inserted hat function $\tilde{\varphi}_{j-1,j}$ and an update of rows $e_{:,j},\:e_{:,j-1}$ and columns $e_{j,:},\:e_{j-1,:}$, highlighted in blue in Figure \ref{fig:spatial_block_ref_ex}.\\
Let us note that this spatial procedure has to be repeated for all the temporal blocks in matrix \eqref{linear_system_explicit}, but nevertheless, as explained above, the modification involves just block elements related to basis with supports containing refined mesh elements, allowing to store and save the remaining of the entries, that can be reused in the next level of refinement.
\begin{figure}
\centering
\includegraphics[scale=0.8]{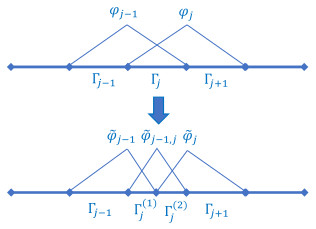}
\caption{Visualization of spatial refinement procedure for piecewise linear basis functions.}
\label{fig:spatial_mesh_ref_ex}
\end{figure}
\begin{figure}
\centering
\includegraphics[scale=0.7]{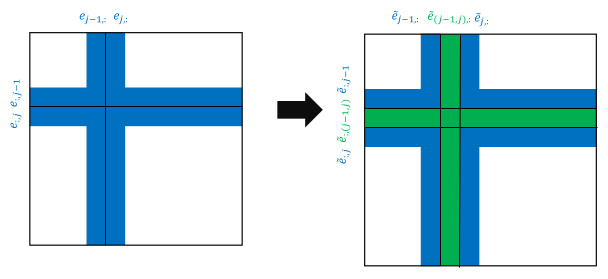}
\caption{Update of a generic block $\mathbb{E}_{\tilde{i},i}$ in \eqref{linear_system_explicit}.}
\label{fig:spatial_block_ref_ex}
\end{figure}

\subsubsection{Wave scattering by a straight crack}\label{Wave scattering by a straight crack}

The first experiment deals with the adaptive resolution of the discrete Neumann problem \eqref{NPdisc}, taking into account the flat obstacle $\Gamma=\left\lbrace x=(x_1,0)\in\mathbb{R}^2\:|\: x_1\in[-0.5,0.5]\right\rbrace$ where the Neumann spatially homogeneous  condition $f(x,t)=H(t)$ is imposed.
The time interval of interest is $[0,T]=[0,2]$. For the following results, the reference energy value $0.79280$ has been obtained by extrapolation from values coming from space-time uniform refinements.\\

At first, let us consider an adaptive approach, where starting from a uniform coarse mesh in space and time with $\Delta t=h=0.1$, we consider an adaptive refinement in space where each marked space element is divided in two equivalent sub-elements, while the time step is uniformly halved at each refinement step. The local coefficient in front of  $\|\mathcal{R}\|_{0,0,  {\Box_{i,j}}}^2$ has been chosen as $\sqrt{\Delta t_i^2+h_j^2}$,
which is equivalent to the one in \eqref{aposteriori} as already remarked, with $s=\frac{1}{2}$.
Results are collected in Figure \ref{energy_error_segment_0} and compared with a space-time uniform refinement approach. The adaptive strategy, as expected, converges faster than the uniform one and, in both cases, the error indicators are capable to perfectly predict the related error decays.\\
\begin{figure}
\centerline{\includegraphics[scale=0.35]{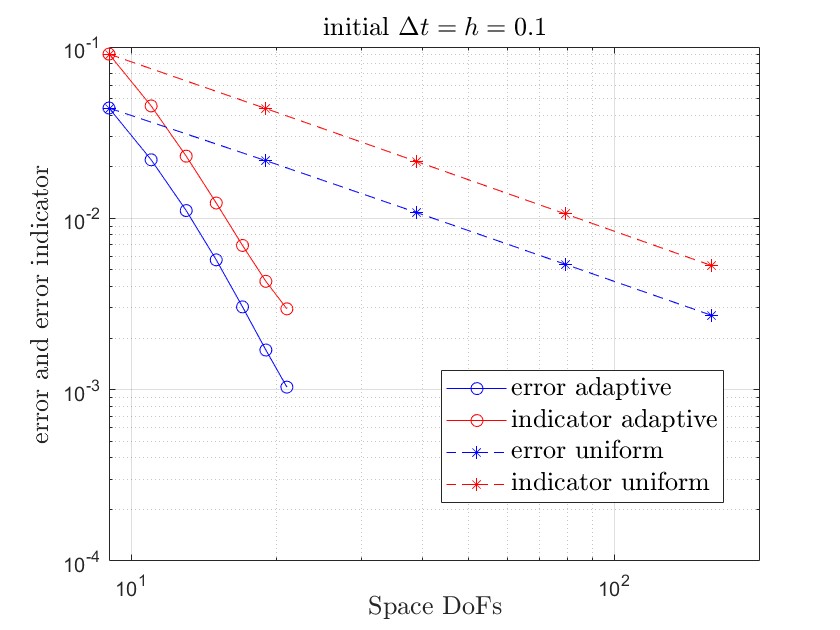}}
\caption{Squared errors in energy norm and related error indicators for the comparison between adaptive and uniform refinement procedures, starting from $\Delta t=h=0.1$. Experiment \ref{Wave scattering by a straight crack}.}
\label{energy_error_segment_0}
\end{figure}

Now let us consider a fixed uniform time discretization with time step $\Delta t=0.0125$, again proceeding with an adaptive space refinement. Discretization in space starts dividing $\Gamma$ uniformly with step $h=0.1$.
Results collected in Figure \ref{energy_error_segment} show a comparison between the squared errors in energy norm given by the adaptive procedure and an iterative uniform refinement of the obstacle: starting from the same space-time level of the discretization, the latter procedure is characterized by a smaller rate of convergence in terms of spatial DoFs.
\begin{figure}
\centerline{\includegraphics[scale=0.35]{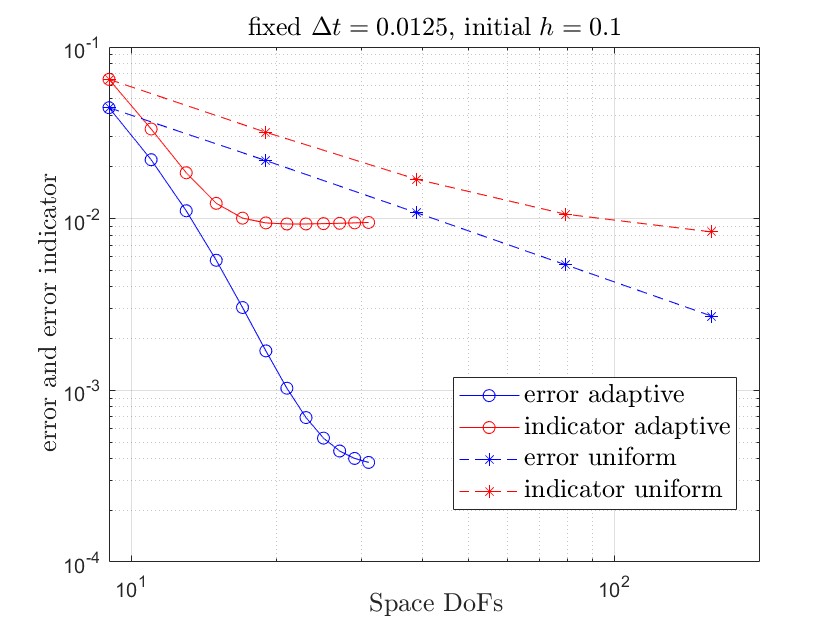}}
\caption{Squared errors in energy norm and related error indicators for the comparison between adaptive and uniform space refinement procedures, having fixed $\Delta t=0.0125$. Experiment \ref{Wave scattering by a straight crack}.}\label{energy_error_segment}
\end{figure}
The adaptive approach produces a fast decay of the energy error, until it stagnates, because the error in time becomes dominant for fixed $\Delta t$ compared to the error reduction from the adaptive space refinements. The error indicator based on the theoretical upper bound \eqref{aposteriori} similarly shows stagnation after an initial decay, both for the space-adaptive and uniform refinement strategies. The stagnation begins approximately when the smallest spatial element of the refinement becomes smaller than the fixed $\Delta t$. 
\\

To complete the analysis, since for this example the dominating error is in space variable, we present numerical results obtained substituting the local coefficient in front of the residuum norm in \eqref{aposteriori} by $h_j$, collected in
Figure \ref{Ale_1_h} for fixed $\Delta t=0.0125$: we note that the $h-$error indicator is capable to strictly follow the behavior of the error.\\
\begin{figure}
\centerline{\includegraphics[scale=0.35]{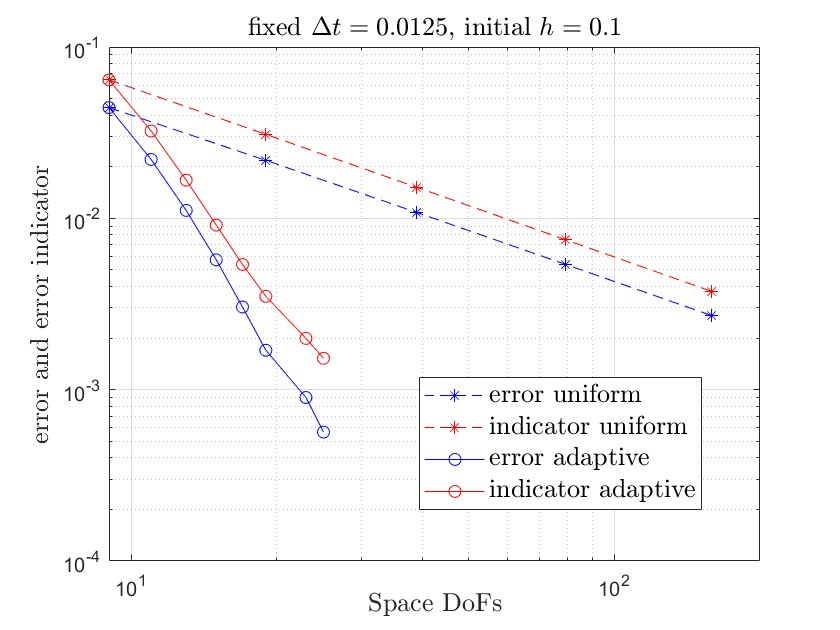}}
\caption{Squared energy errors and $h-$error indicators for fixed  
$\Delta t=0.0125$, for both uniform and adaptive space refinement procedures. Experiment \ref{Wave scattering by a straight crack}.}
\label{Ale_1_h}
\end{figure}
%
 
Finally, in Figure \ref{adaptive refinement slit} it is possible to visualize some of the space meshes generated by the adaptive refinement procedures, where nodes are added near the crack tips of the straight obstacle, which represent the greatest source of error due to the asymptotic behavior of the solution $\psi_{h,\Delta t}$ of problem \eqref{NPdisc}. Moreover, the decreasing values of the $h$-error indicator, related to each mesh element, are shown. 
\begin{figure}
\centerline{\includegraphics[scale=0.45]{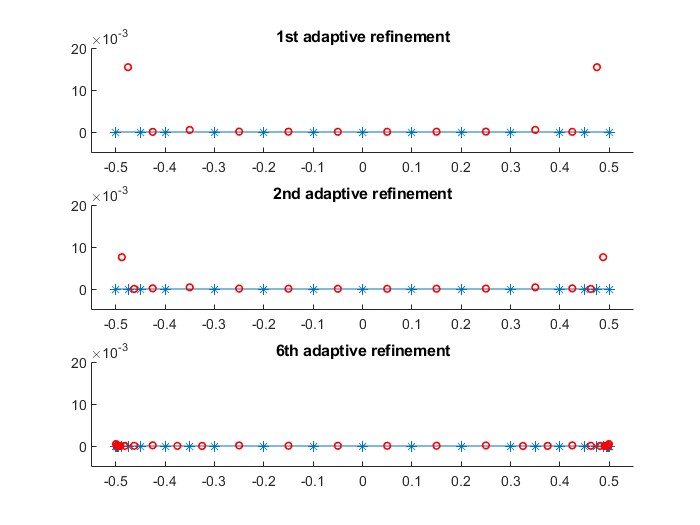}}
\caption{Some space mesh refinement levels obtained using the adaptive procedure (mesh nodes identified by blue asterisks) and $h$-error indicator values (identified by red circles) over mesh elements. Experiment \ref{Wave scattering by a straight crack}.}
\label{adaptive refinement slit}
\end{figure}
 
\subsubsection{Wave scattering by an angular crack}\label{Wave scattering by an angular crack}

The experiment has been repeated considering the open arc $\Gamma= S_1 \cup S_2$, where
$$\begin{array}{c}
S_1=\left\lbrace x=(x_1,x_2)\in\mathbb{R}^2\:\vert \: x_1\in[-0.1,0],\: x_2=(0.1+x)\tan(\pi/3) \right\rbrace\\
     S_2=\left\lbrace x=(x_1,x_2)\in\mathbb{R}^2\:\vert \: x_1\in[0,0.1],\: x_2=(0.1-x)\tan(\pi/3) \right\rbrace  
\end{array}.$$
The angular conjunction of the two segments at point $(0,0.1\tan(\pi/3))$ is an additional source of error since it leads to an asymptotic singular behavior of the BIE solution. The Neumann condition considered is $f(x,t)=n_{x,1}H(t)$ and the time interval of analysis is set as $[0,T]=[0,0.5]$.\\
Starting from a uniform coarse mesh in space and time, fixing at the beginning $\Delta t=h=0.05$, we consider an adaptive refinement in space where each marked space element is divided in two equivalent sub-elements, while the time step is uniformly halved at each refinement step. The local coefficient in front of  $\|\mathcal{R}\|_{0,0,  {\Box_{i,j}}}^2$ has been fixed as $\max\{\Delta t_i,h_j\}$, as given in \eqref{aposteriori}.\\ The reference energy value $0.034012$ has been obtained by extrapolation from values coming from space-time uniform refinements as done in the previous section. Adaptive refinement results are collected in Figure \ref{energy_error_angle} and compared with a space-time uniform refinement approach. The adaptive strategy, as expected, converges faster than the uniform one and, in both cases, the error indicators are capable to perfectly predict the related error decays. Moreover, the expected rate of convergence of $O(\rm{DoFs}^{-1})$ is obtained by uniform refinement, while the adaptive strategy more than doubles the rate of convergence.
In Figure \ref{adaptive refinement angle} some space mesh refinement levels obtained using the adaptive procedure are shown, highlighting, above all at the beginning, the insertion of nodes towards the crack tips and the angle, as expected.
\\
\begin{figure}
\centerline{
\includegraphics[scale=0.35]{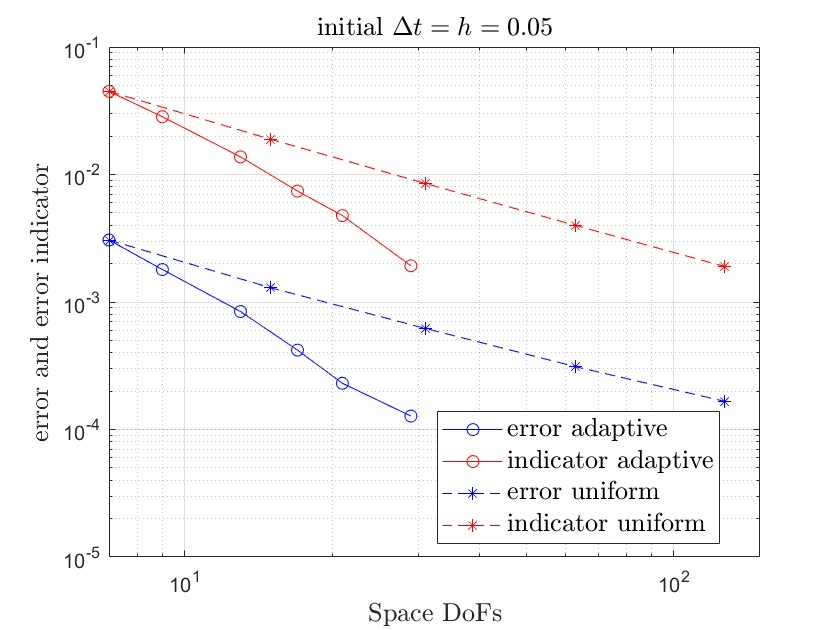}}
\caption{Squared errors in energy norm and related error indicators for the comparison between adaptive and uniform refinement procedures, starting from $\Delta t=h=0.05$. Experiment \ref{Wave scattering by an angular crack}.}
\label{energy_error_angle}
\end{figure}
\begin{figure}[h!]
\hspace{0.1in}
\stackunder[5pt]{\includegraphics[width=0.3\textwidth]{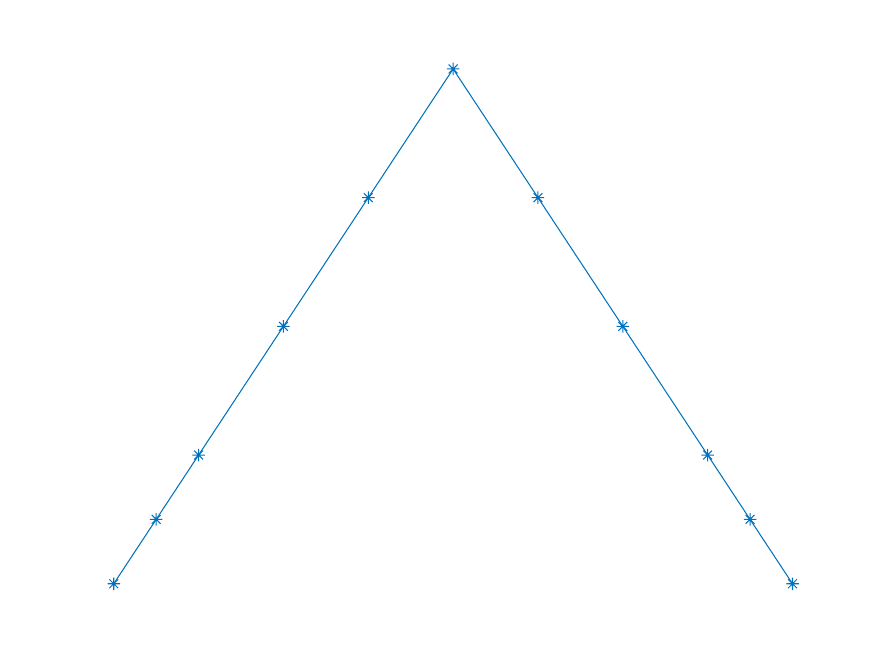}}{1st adaptive refinement}%
\:
\stackunder[5pt]{\includegraphics[width=0.3\textwidth]{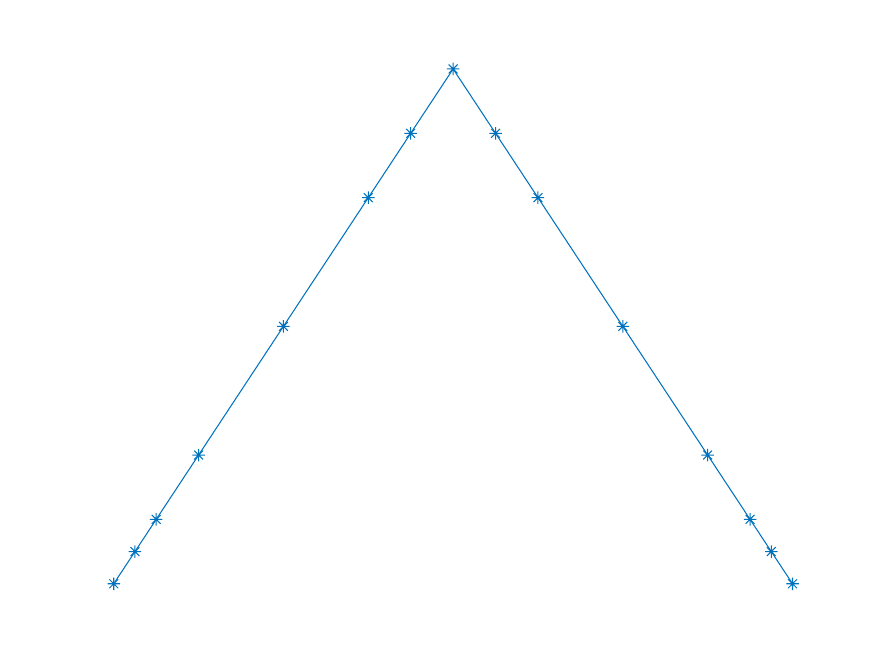}}{2nd adaptive refinement}
\:
\stackunder[5pt]
{\includegraphics[width=0.3\textwidth]{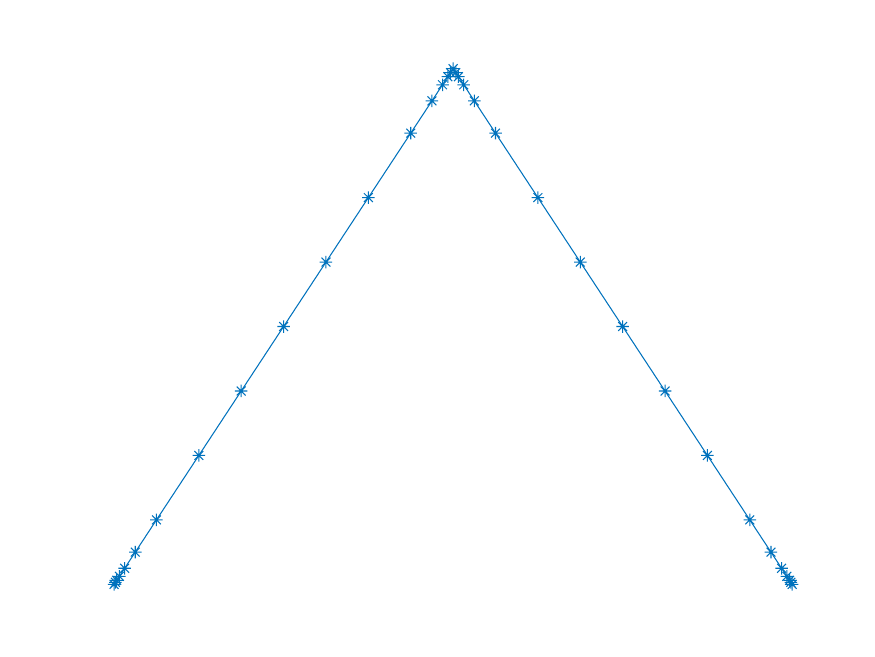}}{6th adaptive refinement}
\vskip 0.5cm
\caption{Some space mesh refinement levels obtained using the adaptive procedure. Experiment \ref{Wave scattering by an angular crack}.}
\label{adaptive refinement angle}
\end{figure}

\subsubsection{Wave scattering by an equilateral triangle}\label{Wave scattering by an equilateral triangle}
\label{Experiment 2: scalar wave scattering by an equilateral triangle}

The second experiment is related to the adaptive resolution of the discrete Neumann problem \eqref{NPdisc} in a domain exterior to the closed equilateral triangle $\Gamma=S_0\cup S_1 \cup S_2$, where
$$\begin{array}{c}
     S_0=\left\lbrace x=(x_1,0)\in\mathbb{R}^2\:\vert \: x_1\in[-0.1,0.1] \right\rbrace,  \\
     S_1=\left\lbrace x=(x_1,x_2)\in\mathbb{R}^2\:\vert \: x_1\in[0,0.1],\: x_2=(0.1-x)\tan(\pi/3) \right\rbrace,  \\ 
     S_2=\left\lbrace x=(x_1,x_2)\in\mathbb{R}^2\:\vert \: x_1\in[-0.1,0],\: x_2=(0.1+x)\tan(\pi/3) \right\rbrace,
\end{array}.$$
The time interval $[0,T]=[0,0.5]$ is considered, and
the following Neumann condition is imposed on $\Gamma$: 
$$
f(x,t)=(\chi_{S_1}(x)-\chi_{S_2}(x))g(t), \qquad\textrm{with}\qquad g(t)=\left\{ \begin{array}{c l}
\sin(4\pi t)^2, & t\in[0,1/8],\\
1, & t>1/8,
\end{array}\right. 
$$
where $\chi_S(x)=1$, if $x\in S$, and $=0$ otherwise. 
Starting from a uniform coarse mesh in space and time, fixing at the beginning $\Delta t=h=0.05$, we consider an adaptive refinement in space where each marked space element is divided in two equivalent sub-elements, while the time step is uniformly halved at each refinement step. The local coefficient in front of  $\|\mathcal{R}\|_{0,0,  {\Box_{i,j}}}^2$ has been fixed as $\max\{\Delta t_i,h_j\}$, as given in \eqref{aposteriori}. The reference  energy value $0.063334$ has been obtained by extrapolation from values coming from space-time uniform refinements, as done in the previous simulations.\\ Adaptive refinement results are collected in Figure \ref{energy_error_eqtr_0} and compared with a space-time uniform refinement approach. The adaptive strategy, as expected, converges faster than the uniform one and, in both cases, the error indicators are capable to perfectly predict the related error decays. Moreover, the expected rate of convergence of  $O(\rm{DoFs}^{-6/5})$ is obtained by uniform refinement. The convergence rate obtained by the adaptive strategy is approximately doubled in this example. In Figure \ref{adaptive refinement triangle_bis} some space mesh refinement levels obtained using the adaptive procedure are shown, highlighting refinements towards the triangle vertices where the given datum presents jumps: the higher the jump, the finer the local mesh.
\\
\begin{figure}
\centerline{\includegraphics[scale=0.35]{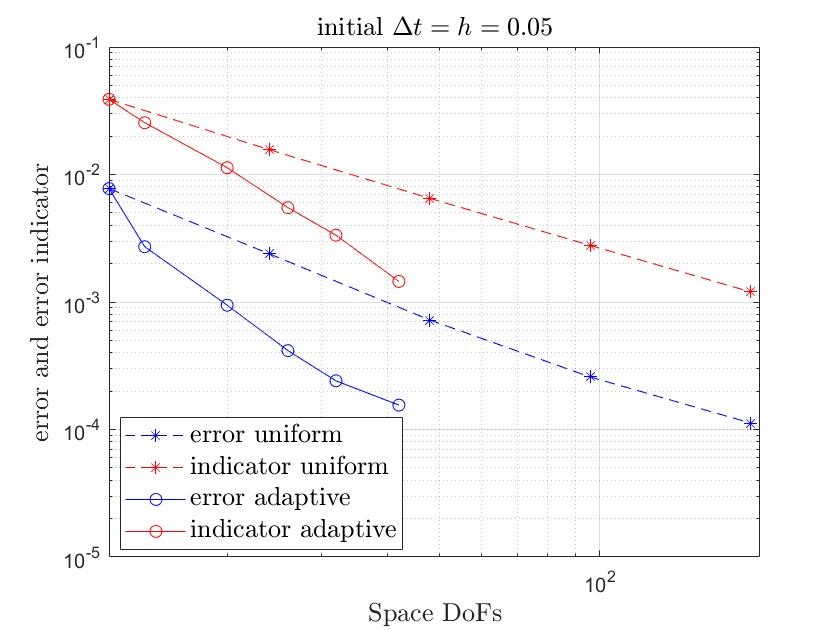}}
\caption{Squared errors in energy norm and related error indicators for the comparison between adaptive and uniform refinement procedures, starting from $\Delta t=h=0.05$. Experiment \ref{Wave scattering by an equilateral triangle}.}
\label{energy_error_eqtr_0}
\end{figure}
\begin{figure}[h!]
\hspace{-0.5in}
\stackunder[5pt]{\includegraphics[width=0.35\textwidth]{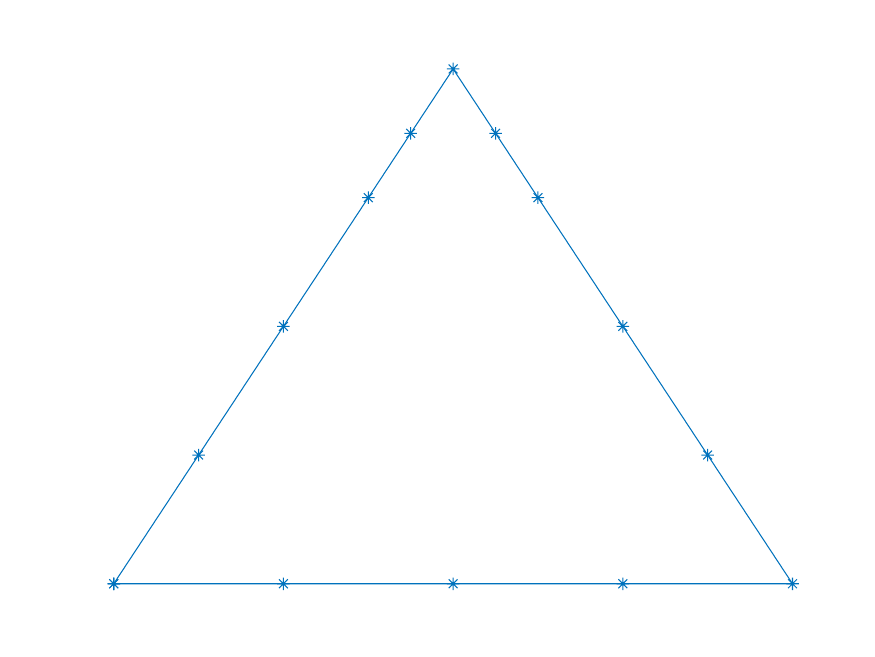}}{1st adaptive refinement}%
\:
\stackunder[5pt]{\includegraphics[width=0.35\textwidth]{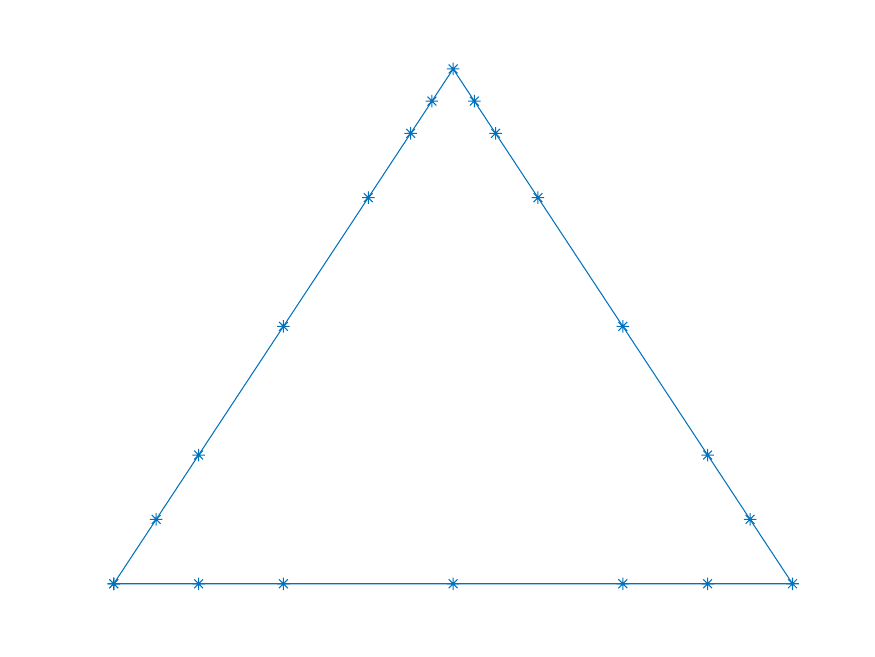}}{2nd adaptive refinement}
\:
\stackunder[5pt]{\includegraphics[width=0.35\textwidth]{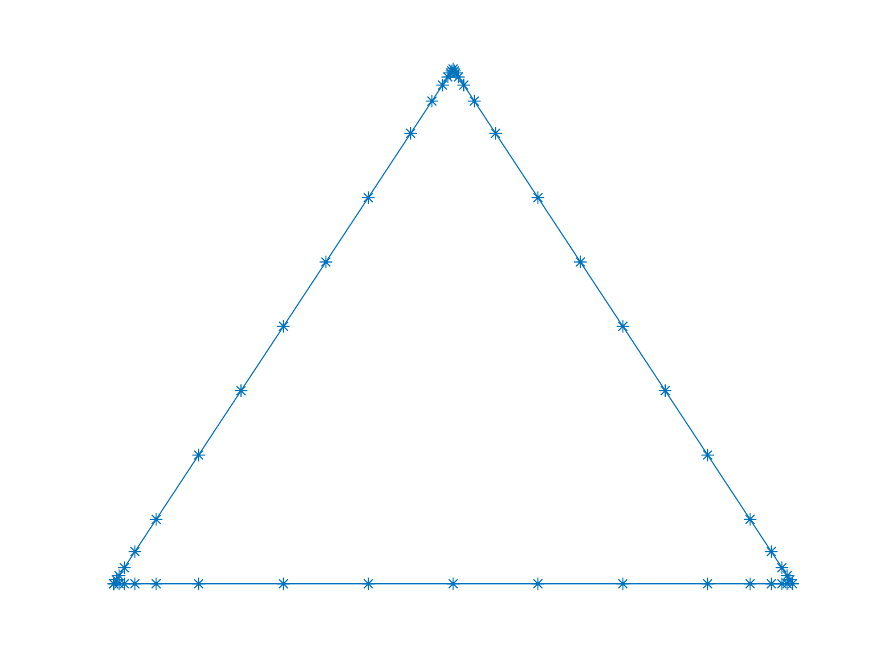}}{6th adaptive refinement}
\vskip 0.5cm
\caption{Some space refinements levels obtained using the adaptive procedure.  Experiment \ref{Wave scattering by an equilateral triangle}.}
\label{adaptive refinement triangle_bis}
\end{figure}

\subsection{Time-adaptivity}\label{subsection:timeadaptive}

Non-uniform choices of the time step are relevant to approximate singular behavior of the solution in time. Prescribed, graded time step sizes were first considered in \cite{ss}. Adaptive refinements of the time step size were recently studied for the Dirichlet problem in one space dimension \cite{hoonhout2023,zank1d}, after earlier work in \cite{sv}. Analogously, as a second step towards adaptive mesh refinements for the hypersingular integral equation in space and time, we consider  refinements of the temporal mesh, based on space-integrated error indicators. The precise algorithm is given in the following:\\

\noindent {{\textbf{Time-Adaptive Algorithm}
}\\
 Input: Mesh $\mathcal{T}=(\mathcal{T}_T \times \mathcal{T}_\Gamma)_0$, refinement parameter $\theta  \in (0, 1)$, tolerance $\epsilon > 0$, datum $f$.
\begin{enumerate}
\item Solve \eqref{NPdisc} on $\mathcal{T}$.
\item Compute the error indicators $\eta(\Box_{i,j})$  
in each space-time element $\Box_{i,j} \in \mathcal{T}$, i.e., for $i=0, \dots, N_T-1$ and $j=1,\dots, M_\Gamma$. 
\item Evaluate $\eta_i=\displaystyle \sum_{j=1}^{M_\Gamma} \eta(\Box_{i,j}), \, i=0,\dots,N_T-1$.
\item Stop if $\displaystyle \sum_{i=0,\dots,N_T-1} \eta_i < \epsilon$\,, otherwise:
\item Find $\displaystyle \eta_{max}=\max_{i=0,\dots,N_T-1} \,\eta_i$.
\item Mark all $\Box_{i,j}\in \mathcal{T}$ with $i$ such that $\eta_i > \theta \,\eta_{max}$ and $j=1, \dots, M_\Gamma$
\item Refine each marked $\Box_{i,j}$ in time, keeping $ \min_i \frac{\Delta t_i}{h_j}$ fixed.
\item Go to 1.
\end{enumerate}}

Note that the time-adaptive algorithm may lead to large CFL ratios $\frac{\Delta t_i}{h_j}$. The Galerkin discretization remains stable in this case. Keeping the minimum of $\frac{\Delta t_i}{h_j}$ fixed assures the reduction of the space discretization error. As for the space-adaptive algorithm, also other choices for the marking are commonly used, e.g.~D\"{o}rfler marking.
Instead of refining in space in Step 7., one can also consider time-adaptive algorithms for a sufficiently small, fixed spatial mesh size.\\ 
For the following numerical simulations, the parameter $\theta=0.4$ has been considered.

\subsubsection{Time refinement step}
Time refinement procedure could induce the splitting of a time interval $I_i$ in the subintervals $I_{i}^{(1)},I_{i}^{(2)}$ and this will produce a substitution of the basis function $\nu_i$, considered in Section \ref{algo_sub},
with the new ramps $\nu_{i}^{(1)},\nu_i^{(2)}$. Let us note that, even if a new time knots has been added inside the support of the adjacent basis function $\nu_{i-1}$, this latter is not affected at all.
\begin{figure}
\centering
\includegraphics[scale=0.8]{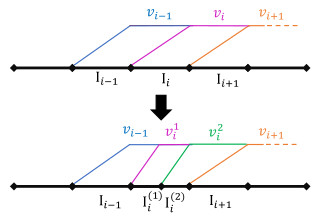}
\caption{Visualization of time refinement procedure for piecewise linear basis functions. }
\label{fig:mesh_time_refinement}
\end{figure}

In \eqref{matrix time ref} we show how time refinement modifies the block structure of the matrix $\mathbb{E}$ in \eqref{linear_system_explicit}: two additional rows and columns blocks related to $\nu_{i}^{(1)},\nu_i^{(2)}$ are added in substitution of the row and column blocks related to $\nu_i$, 
keeping anyway the global lower triangular block structure 
\begin{equation}\label{matrix time ref}
\begin{pmatrix}
\mathbb{E}_{0,0} & \cdots  & 0 &0 & \ldots  & 0\\
\vdots & \ddots  & 0  & 0  & \ldots & 0\\
\mathbb{E}^{(1)}_{i,0} & \cdots & \mathbb{E}^{(1,1)}_{i,i} & 0 &\ldots & 0\\
\mathbb{E}^{(2)}_{i,0} & \cdots & \mathbb{E}^{(2,1)}_{i,i} & \mathbb{E}^{(2,2)}_{i,i} & \ldots & 0\\
\vdots & \ddots & \vdots  & \vdots & \ddots & \vdots\\
\mathbb{E}_{N_T-1,0} & \cdots & \mathbb{E}^{(1)}_{N_T-1,i} & \mathbb{E}^{(2)}_{N_T-1,i} & \ldots & \mathbb{E}_{N_T-1,N_T-1} 
\end{pmatrix}.
\end{equation}
Unlike the refinements of the spatial mesh in Subsection \ref{subsection:spaceadaptive}, the refinements in time lead to a system \eqref{matrix time ref} without any Toeplitz structure.

\subsubsection{Wave scattering by a circle}\label{Wave scattering by a circle}

The last experiment is conceived to properly add temporal mesh nodes. We consider, as boundary geometry of the problem, the smooth closed arc  $\Gamma=\left\lbrace x=(x_1,x_2)\in\mathbb{R}^2\:\vert\:\sqrt{x_1^2+x_2^2}=0.5 \right\rbrace$. Since the considered $\Gamma$ is a circle, the solution of the weak problem \eqref{NPdisc} is smooth in space and, unlike the previous examples, is not dominated by the spatial discretization error. We fix a uniform spatial mesh with $h=\pi/32\simeq 0.1$.\\
The Neumann datum $f(x,t)=(5t)^{-0.27}e^{-5t}$ is imposed: the singularity at $t=0$ will produce an asymptotic behavior of type $\mathcal{O}(t^{\alpha})$ with $0<\alpha<1$ in the solution of \eqref{NPdisc}, so that adaptive mesh refinements in time are expected near $t=0$. The time interval of interest $[0,T]=[0,\pi/4]$ is discretized both uniformly and adaptively, starting with a uniform decomposition of $[0,T]$ made by $8$ intervals, i.e., fixing initial $\Delta t=\pi/32\simeq 0.1$. The reference  energy value $1.777$ has been obtained by extrapolation from energy values computed from uniform refinements in space and time.\\
In Figure \ref{energy_error_time_max} we report {the behavior of the squared energy error and the error indicator, based on \eqref{aposteriori} with $s=\frac{1}{2}$}. We refine the time mesh, both uniformly and adaptively: starting from the same space-time level of the discretization, the adaptive time refinement shows a higher rate of convergence than the time uniform procedure in terms of temporal DoFs. Moreover, even if the local coefficient multiplying the residuum norm in \eqref{aposteriori} remains constant, i.e., $\max\{\Delta t_i,h_j\}=h$, in this case we don't observe a stagnation of the error indicator, as it happens in the example on the straight crack, see Figure \ref{energy_error_segment}. 
In fact, compared to the experiment on the segment, the solution on the circle is smooth in space. 
Correspondingly,  the error indicator decays in the current example, unlike on the segment, even if at a slower rate than the energy error.\\
In Figure \ref{time refinement Ale} we show some time mesh refinement levels obtained using the adaptive procedure and the error indicator values over meshes elements: in the first refinement steps only one time node is added always near $t=0$, as expected.\\
\begin{figure}[h!]
\centering
\includegraphics[width=0.65
\textwidth]
{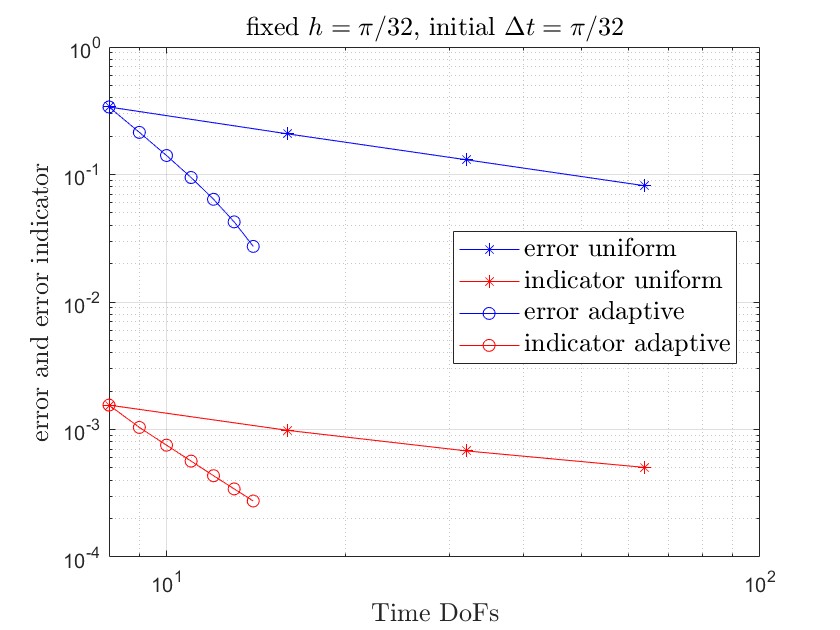}
\caption{Squared energy errors and error indicators for fixed $h=\pi/32$, for both uniform and adaptive time refinements procedures.  Experiment \ref{Wave scattering by a circle}.}
\label{energy_error_time_max}
\end{figure}
\begin{figure}[h!]
\centering
\includegraphics[scale=0.45]{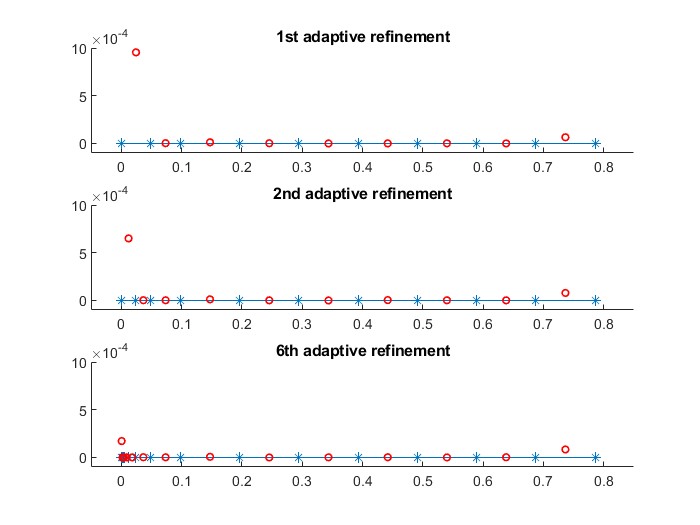}
\caption{Some time mesh refinement levels obtained using the adaptive procedure 
(mesh nodes identified by blue asterisks) and error indicator values (identified by red circles) over mesh elements. Experiment \ref{Wave scattering by a circle}.}
\label{time refinement Ale}
\end{figure}

Since in this example the dominating error is in the time variable, we present also numerical results obtained substituting the local coefficient in front of the norm of the residual in \eqref{aposteriori}, by the heuristic choice $\Delta t_i$. Results of this experiment are collected in Figure \ref{energy_error_time_dt}: starting from the same space-time level of the discretization, the uniform time refinement presents a lower rate of convergence than the time-adaptive procedure in terms of temporal DoFs. However, {differently from the heuristic $h-$indicator of the experiment \ref{Wave scattering by a straight crack},} the $\Delta t-$error indicator decays faster than the error in energy norm, showing that they are not in accordance, even if this heuristic indicator leads to sensible mesh refinements. 
In Figure \ref{time refinement} the temporal nodes obtained by the successive iterations of this latter time-adaptive approach are displayed: as expected, the asymptotic behavior of the solution in $t=0$ leads to an accumulation of the mesh nodes mostly at the beginning of the considered time interval.

\begin{figure}[h!]
\centering
\includegraphics[width=0.6
\textwidth]
{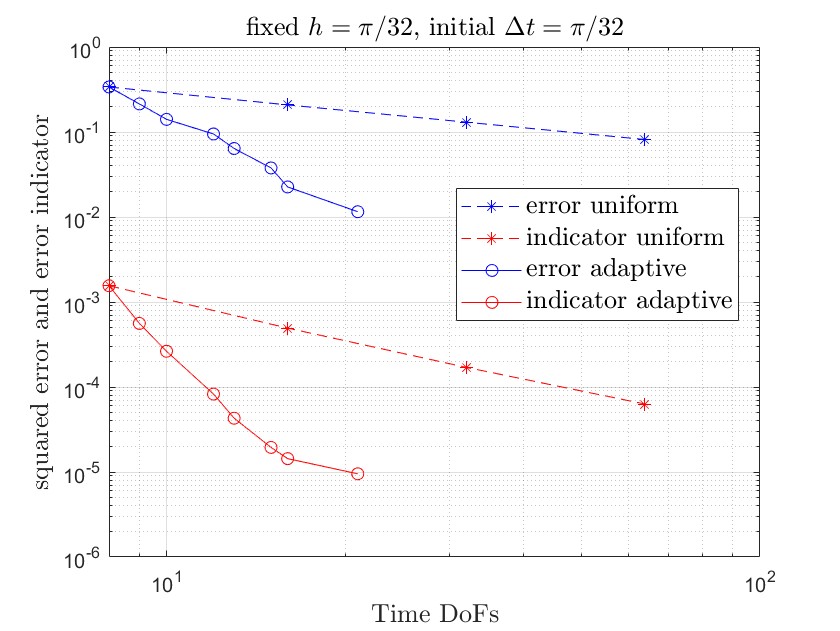}
\caption{Squared energy errors and $\Delta t-$error indicators for fixed $h=\pi/32$, for both uniform and adaptive time refinements procedures. Experiment \ref{Wave scattering by a circle}.}
\label{energy_error_time_dt}
\end{figure}
\begin{figure}[h!]
\centering
\includegraphics[width=1\textwidth]{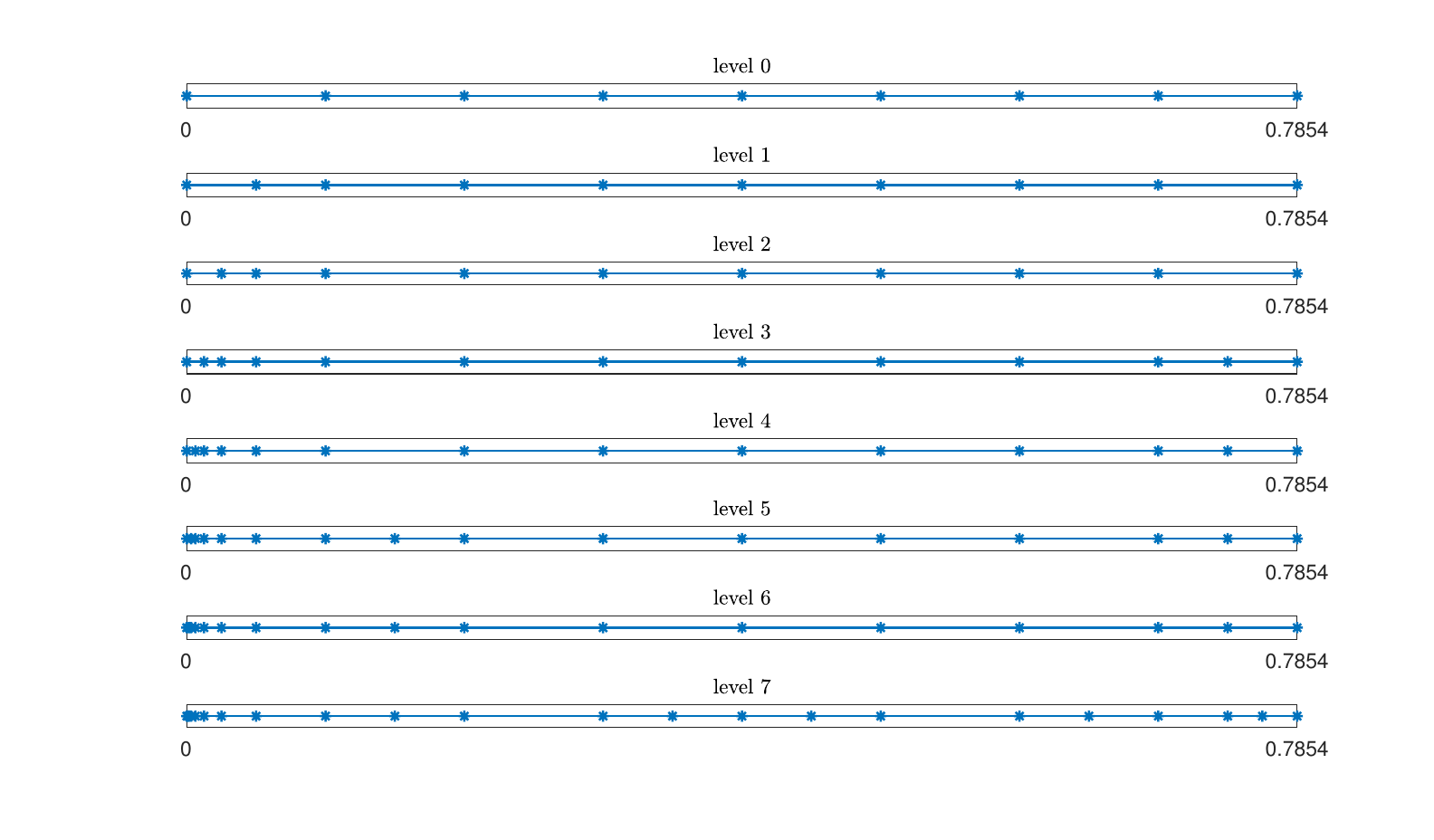}
\caption{Successive refinements of time mesh in the adaptive procedure with $\Delta t-$error indicator. Experiment \ref{Wave scattering by a circle}.}
\label{time refinement}
\end{figure}

\subsection{Memory savings and CPU times of the adaptive procedures}
We use the following quantity to compare the memory requirements of the adaptive and non-adaptive simulations:
\begin{equation}\label{memory saving}
    \mathcal{S}=1-\dfrac{\mathcal{M}_a}{\mathcal{M}_u}\,,
\end{equation}
where $\mathcal{M}_u=M_{\Gamma,u}^2N_{T,u}$ denotes the memory needed to store the Toeplitz matrix in the case of a uniform space-time discretization. The definition of the memory requirement $\mathcal{M}_a$ in the adaptive approach depends on whether time or space refinement are used:
\begin{equation}
    \mathcal{M}_a=\left\lbrace
    \begin{array}{lr}
         M_{\Gamma,a}^2N_{T,a}\,,& \textrm{for space-adaptive procedures,}\\[2pt]
         M_{\Gamma,a}^2\frac{N_{T,a}(N_{T,a}+1)}{2}\,,& \textrm{for time-adaptive procedures.}
    \end{array}
    \right.
\end{equation}
For the space-adaptive simulations in Subsection \ref{Wave scattering by a straight crack}, fixing an error level of $3\cdot10^{-3}$, we find memory savings corresponding to $\mathcal{S}=99\%$. Similarly, for the space-adaptive simulations in Subsections \ref{Wave scattering by an angular crack} and \ref{Wave scattering by an equilateral triangle}, considering an error level of $1\cdot10^{-4}$, we obtain $\mathcal{S}=90\%$. Finally, in the case of the two time-adaptive experiments in Subsection \ref{Wave scattering by a circle}, with fixed error levels of $2\cdot10^{-2}$, respectively $1\cdot10^{-2}$, the memory saving corresponds to $\mathcal{S}=60\%$, respectively $\mathcal{S}=77\%$. Note that we observe memory savings even though the Toeplitz structure of the Galerkin matrix is lost for the time-adaptive procedure. In any case, the more the required accuracy, the higher the memory savings of the adaptive approach with respect to the uniform one.\\
The CPU times of the adaptive algorithm are specific to the implementation, as they strictly depend on code optimization and on the used hardware. While we cannot provide general formulas, for the reader's convenience we refer to two typical simulations on a laptop with Intel(R) Core(TM) i5-10210U CPU, parallelized over 4 cores at 1.60GHz frequency, with 16 GB RAM.\\
For the space-adaptive simulation related to Figure \ref{energy_error_segment}, to reach an error level of $10^{-3}$ the run times were: $577\,s$ for the construction of the linear systems at the involved refinement steps,  $57\,s$ for their solution, $389\, s$ for the evaluation of error indicators and $1\, s$ for the generation of the adaptive meshes.\\
For the time-adaptive simulation related to Figure \ref{energy_error_time_max}, to reach an error level of $10^{-1}$ the run times were: $380\,s$ for the construction of the linear systems at the involved refinement steps, $60\,s$ for their solution, $25\, s$ for the evaluation of error indicators and $1\, s$ for the generation of the adaptive meshes.

\section{Conclusions}

This article presented an a posteriori error estimate for the  time-domain hypersingular integral equation corresponding to the wave equation and studied the resulting adaptive mesh refinement procedures, either in space or in time. As recently explored for Dirichlet boundary conditions, {or the weakly singular integral equation,} \cite{cf23,chaumont2024damped,adaptive} local refinements of the spatial mesh lead to efficient approximations for solutions with time-independent, geometric singularities, while refinements of the time mesh efficiently approximated singular transient behavior. The adaptive procedures, in space, resp.~in time, were studied by extensive numerical simulations for tensor-product meshes.\\
In general, solutions to wave equations naturally exhibit singularities which move in space and time, such as sharp wave crests. Fully space-time-adaptive mesh refinements beyond tensor-product meshes are required to resolve such singularities and are being studied for the weakly singular operator \cite{vdraft,Glaefke}. 
The space-time  adaptive procedure is again based on the four steps:}

\vspace*{-0.5cm}

 \begin{align*}
  \textbf{SOLVE}&\longrightarrow \textbf{ESTIMATE}\longrightarrow \textbf{MARK}\longrightarrow \textbf{REFINE},
\end{align*}

\noindent where a posteriori error estimates are used to define error indicators in each element of the space-time mesh. Hence, the flexibility of the space-time adaptive approach comes with additional computational
cost. Moreover, in this general scenario, the Galerkin matrix no longer has a lower triangular block structure, so that the traditional time-stepping approach has to be replaced by the solution of the full space-time system.\\
The special block Toeplitz structure of the Galerkin matrix is, naturally, also lost when the time step is not uniform or the spatial mesh change with time, increasing the memory requirement, even if this latter remains less than the memory requirements of a uniform refinement approach, considering the same level of accuracy of the approximate solution.\\ Efficient
implementations of fully space-time adaptive boundary elements should reuse those matrix entries which have not been affected by a refinement step, but the actual efficient implementation remains a challenge of ongoing work \cite{vdraft}.

For the hypersingular operator, the a posteriori error estimate in this article defines residual error indicators $\eta(\Box_{i,j})$ in each space-time element on general meshes, which allows to formulate fully space-time adaptive algorithms. Unlike for the weakly singular operator, the accurate quadrature of the hypersingular integrals in the assembly of the Galerkin matrix, however, differs significantly from existing implementations on tensor product meshes and will be addressed in future research.

\vspace*{-0.3cm}

\section*{Acknowledgments}

This research was supported by the  Research in Residence program of the Centre International de Rencontres Math\'{e}matiques, Luminy, in 2023.

\vspace*{-0.3cm}

\section{Appendix}

Space-time anisotropic Sobolev spaces provide a convenient scale of function spaces to study  layer potentials  on the boundary $\Gamma$. We refer to \cite{hd, setup} for an extended exposition. To define these spaces on a screen, i.e., if $\partial\Gamma\neq \emptyset$, we first extend $\Gamma$ to a closed, orientable Lipschitz manifold $\widetilde{\Gamma}$. 

On $\Gamma$ one we consider the usual Sobolev spaces of supported distributions:

\vspace*{-0.2cm}

$$\widetilde{H}^s(\Gamma) = \{u\in H^s(\widetilde{\Gamma}): \mathrm{supp}\ u \subset {\overline{\Gamma}}\}\ , \quad\ s \in \mathbb{R}\ .$$
The space ${H}^s(\Gamma)$ is then the quotient $ H^s(\widetilde{\Gamma}) / \widetilde{H}^s({\widetilde{\Gamma}\setminus\overline{\Gamma}})$.} \\
We now write down an explicit family of Sobolev norms. To do so, fix a partition of unity $\alpha_i$ subordinate to a covering of $\widetilde{\Gamma}$ by open sets $B_i$. For diffeomorphisms $\varphi_i$ mapping each $B_i$ into the unit cube $\subset \mathbb{R}^d$, a family of norms with parameter $\omega \in \mathbb{C}\setminus \{0\}$ is given by:
\begin{equation*}
 \|u\|_{s,\omega,{\widetilde{\Gamma}}}=\left( \sum_{i=1}^p \int_{\mathbb{R}^n} (|\omega|^2+|\xi|^2)^s|\mathcal{F}\left\{(\alpha_i u)\circ \varphi_i^{-1}\right\}(\xi)|^2 d\xi \right)^{\frac{1}{2}}\ .
\end{equation*}
Here, $\mathcal{F}$ denotes the Fourier transform. The norms for different parameters $\omega \in \mathbb{C}\setminus \{0\}$ are equivalent. 

The above norms induce Sobolev norms on the spaces $H^s(\Gamma)$, $\|u\|_{s,\omega,\Gamma} = \inf_{v \in \widetilde{H}^s(\widetilde{\Gamma}\setminus\overline{\Gamma})} \ \|u+v\|_{s,\omega,\widetilde{\Gamma}}$ and $\widetilde{H}^s(\Gamma)$, $\|u\|_{s,\omega,\Gamma, \ast } = \|e_+ u\|_{s,\omega,\widetilde{\Gamma}}$. Here, $e_+$ extends the distribution $u$ by $0$ from $\Gamma$ to $\widetilde{\Gamma}$.  Note that the norm $\|u\|_{s,\omega,\Gamma, \ast }$ corresponds to extension by zero, while $\|u\|_{s,\omega,\Gamma}$ allows  extension by an arbitrary $v$. Therefore, $\|u\|_{s,\omega,\Gamma, \ast }$ is stronger than $\|u\|_{s,\omega,\Gamma}$. 

When $H^s(\Gamma)$ is endowed with  the norm for a fixed parameter $\omega$, we denote it by $H^s_\omega(\Gamma)$. Similarly, $\widetilde{H}^s_\omega(\Gamma)$ denotes the space $\widetilde{H}^s(\Gamma)$ with the norm corresponding to $\omega$. 

We can now define the class of space-time anisotropic Sobolev spaces relevant for the main body of this article:
\begin{definition}
For $\sigma>0$ and $r,s \in\mathbb{R}$ define
\begin{align*}
 H^r_\sigma(\mathbb{R}^+,H^s(\Gamma))&=\{ u \in \mathcal{D}^{'}_{+}(H^s(\Gamma)): e^{-\sigma t} u \in \mathcal{S}^{'}_{+}(H^s(\Gamma))  \textrm{ and }   ||u||_{r,s,\sigma,\Gamma} < \infty \}\ , \\
 H^r_\sigma(\mathbb{R}^+,\widetilde{H}^s({\Gamma}))&=\{ u \in \mathcal{D}^{'}_{+}(\widetilde{H}^s({\Gamma})): e^{-\sigma t} u \in \mathcal{S}^{'}_{+}(\widetilde{H}^s({\Gamma}))  \textrm{ and }   ||u||_{r,s,\sigma,\Gamma, \ast} < \infty \}\ .
\end{align*}
Here, $\mathcal{D}^{'}_{+}(E)$ respectively~$\mathcal{S}^{'}_{+}(E)$ denote the spaces of distributions, respectively tempered distributions, on $\mathbb{R}$ with support in $[0,\infty)$, taking values in $E$, where $E= {H}^s({\Gamma})$, respectively  $E=\widetilde{H}^s({\Gamma})$. The space-time anisotropic Sobolev norms are given by
\begin{align*}
\|u\|_{r,s,\sigma}:=\|u\|_{r,s,\sigma,\Gamma}&=\left(\int_{-\infty+i\sigma}^{+\infty+i\sigma}|\omega|^{2r}\ \|\hat{u}(\omega)\|^2_{s,\omega,\Gamma}\ d\omega \right)^{\frac{1}{2}}\ ,\\
\|u\|_{r,s,\sigma,\ast} := \|u\|_{r,s,\sigma,\Gamma,\ast}&=\left(\int_{-\infty+i\sigma}^{+\infty+i\sigma}|\omega|^{2r}\ \|\hat{u}(\omega)\|^2_{s,\omega,\Gamma,\ast}\ d\omega \right)^{\frac{1}{2}}\,.
\end{align*}
\end{definition}
When $r=s=0$, $H^0_\sigma(\mathbb{R}^+,H^0(\Gamma))$ and $H^0_\sigma(\mathbb{R}^+,\widetilde{H}^0(\Gamma))$ both coincide with the weighted $L^2$-space defined by the scalar product $\langle u,v \rangle_\sigma := \int_0^\infty e^{-2\sigma t} \int_\Gamma u \overline{v} ds_x\ dt$.

\vspace*{-0.4cm}

\bibliographystyle{plain}
\bibliography{Biblio}

\end{document}